\documentclass{report}

\usepackage{amsfonts,amsmath,eufrak}

\def\frak{\EuFrak} 
\def\nid{\noindent}
\def\Hom{\mathop{\rm Hom}\nolimits}
\def\Der{\mathop{\rm Der}\nolimits}
\def\End{\mathop{\rm End}\nolimits}
\def\Ker{\mathop{\rm Ker}\nolimits} 
\def\Ext{\mathop{\rm Ext}\nolimits}
\def\Deg{\mathop{\rm Deg}\nolimits}
\def\Tor{\mathop{\rm Tor}\nolimits}
\def\Gr{\mathop{\rm Gr}\nolimits}
\def\sgn{\mathop{\rm sgn}\nolimits}
\def\ad{\mathop{\rm ad}\nolimits}
\def\ds{\displaystyle}
\def\C{\mathcal {C}}
\def\D{\mathcal {D}}
\def\qed{{\hfill\vrule height 5pt width 5pt depth 0pt}}
\def\phi{\varphi}
\def\Im{\mathop {\rm Im}}
\def\induced{\mathop{\rm Ind}\nolimits}

\newcommand{\bino}[2]{\ensuremath{\left ({#1}\atop{#2}\right )}}

\newtheorem{lemma}{Lemma}[chapter]
\newtheorem{theorem}[lemma]{Theorem}
\newtheorem{definition}[lemma]{Definition}
\newtheorem{corollary}[lemma]{Corollary}
\begin{document}
\pagenumbering{roman}
\begin{center}
{\Large \bf Cohomology of Restricted Lie Algebras}
\vskip 0.2in
\centerline{By}
\vskip 0.1in
\centerline{\large  TYLER JONAH EVANS}
\centerline{\large B.A. (Sonoma State University) 1994}
\centerline{\large M.S. (University of Oregon) 1996}
\vskip 0.3in
\centerline{\large  DISSERTATION}
\vskip 0.1in
\centerline{\large Submitted in partial satisfaction of the requirements for the
degree of}
\vskip 0.1in
\centerline{\large  DOCTOR OF PHILOSOPHY}
\vskip 0.1in
\centerline{\large in}
\vskip 0.1in
\centerline{\large MATHEMATICS}
\vskip 0.2in
\centerline{\large in the}
\vskip 0.1in
\centerline{\large  OFFICE OF GRADUATE STUDIES}
\vskip 0.1in
\centerline{\large of the}
\vskip 0.1in
\centerline{\large  UNIVERSITY OF CALIFORNIA}
\vskip 0.1in
\centerline{\large  DAVIS}
\end{center}
\vskip 0.1in
{\large Approved:}
\begin{center}
\centerline{\underbar{\hskip 2.5in}}
\vskip 0.15in
\centerline{\underbar{\hskip 2.5in}}
\vskip 0.15in
\centerline{\underbar{\hskip 2.5in}}
\vskip 0.2in
\centerline{\large Committee in Charge}
\vskip 0.2in
\centerline{\large 2000}
\end{center}
\newpage
\large
\tableofcontents


\newpage
\large
{\Large \bf ACKNOWLEDGEMENTS} \\

I would like to thank my advisor and friend Professor Dmitry B. Fuchs, without whom this dissertation could not have been possible.  His seemingly infinite patience, especially during this last year, was invaluable in completing this work.  I am forever indebted to him for the hours of stimulating conversations in which he guided this research.  His passion for, and insight into, the beauty of mathematics is truly inspiring, and I am honored to have had the opportunity to learn mathematics from him.

Many of the professors at UC Davis contributed to my education, and I am grateful to each of them. I would like to especially thank 
Professor Albert Schwarz for the outstanding lecture courses and illuminating 
 conversations that I enjoyed during my study at UC Davis. I would also like to thank Professor Motohico Mulase for his kind attention, both personal and mathematical, during my graduate study at UC Davis.

I would like to thank the entire staff in the Department of Mathematics for their support over the past four years, especially Mrs.\  Kathy La Giusa whose assistance  during my years of graduate study was simply indispensable. 

I would like to express my appreciation to Dr.\  Michael Penkava for many useful conversations and suggestions, and I am grateful to Mr.\  Peter Littig for reading and (gently) correcting the manuscript.

I am grateful to all of my friends and colleagues in Davis.  Their unswerving support always encouraged me to hold onto my dreams.  

I would like to thank my parents and my brother, Mr.\  M. Bryan Evans, whose support over the years has enabled me to achieve this goal.  I am truly grateful for their love and support.  Finally, I want to thank my life long friend Ms.\  Holly Evans; what a long strange trip it's been!


\newpage 
\begin{center}
\underline{\bf \Large Abstract} 
\end{center}

\medskip

\large
\noindent 
In this dissertation, we investigate the cohomology theory of restricted Lie algebras.  Motivations for the definition of a restricted Lie algebra are given and the theory of ordinary Lie algebra cohomology is briefly reviewed, including a discussion on algebraic interpretations of the low dimensional cohomology spaces of ordinary Lie algebras.  The general Cartan-Eilenberg construction of the standard cochain complex is given for ordinary Lie algebras.  The representation theory of restricted Lie algebras is reviewed including a description of the restricted universal enveloping algebra $U_{\rm res.}(\frak g)$ of a restricted Lie algebra.  In the case of an abelian restricted Lie algebra, we construct an augmented complex of free $U_{\rm res.}(\frak g)$ modules that is exact in dimensions less than $p$ and hence define the cohomology theory of these algebras in dimension less than $p$.  Explicit formulas for the dimensions of the cochain spaces are given in the abelian case. In particular we show that the dimension of $C^k(\frak g)$ is the same as that of the symmetric algebra $S^k(\frak g)$.  In the non-abelian case, we explicitly construct a cochain complex $\{C^k(\frak g;M),\delta^k\}$ for any coefficient module $M$ for $k\le 3$ and give explicit formulas for the coboundary operators in these dimensions.  It is shown that classical and restricted cohomology do not differ at all in dimension zero and that the restricted cohomology space $H^1(\frak g;M)$ is canonically injected into  the classical cohomology $H^1_{\rm cl.}(\frak g;M)$.  A canonical map $H^2(\frak g;M)\to H^2_{\rm cl.}(\frak g;M)$ is constructed and the kernel is investigated for specific coefficient modules.  The corresponding notions of the usual algebraic interpretations of ordinary low dimensional cohomology are defined and we show that our restricted cohomology spaces encode this information as well. The dissertation concludes with some remarks on multiplicative structures in our complex as well as directions for further research.


\newpage
\pagestyle{myheadings} 
\pagenumbering{arabic}
\markright{  \rm \normalsize CHAPTER 1. \hspace{0.5cm}
 COHOMOLOGY OF RESTRICTED LIE ALGEBRAS }
\large 
\chapter{Introduction}
\thispagestyle{myheadings}
The theory of Lie groups and their Lie algebras was first developed by Sophus Lie in the latter part of the nineteenth century as a part of geometry.  During the period from 1900 to 1940, more and more of the theory of Lie algebras of characteristic zero was developed purely algebraically.  Weyl's Theorem (1925) was originally proved using integration on compact groups, and it was just ten years later when Casmir and van der Waerden found a purely algebraic proof.  This work together with J.H.C. Whitehead's two lemmas were among the hints that enabled Chevalley and Eilenberg to construct the cohomology spaces $H^*(\frak g;M)$.

In the positive characteristic case, many  Lie algebras that arise as natural examples possess and additional structure. The guiding example is the algebra $\Der(A)$ of derivations of an associative algebra $A$.  It is well known that the Lie commutator $[D,D']$ of two derivations is again a derivation, but the composition $D\circ D'$ need not be a derivation.  However, if $D$ is a derivation of $A$ and $k\ge 1$, we do have the Leibniz rule
\[
D^k(ab)=\sum_{j=0}^k \bino {k}{j}D^j(a)D^{k-j}(b).
\]
In particular, if the characteristic $p$ of the ground field is positive, and we take $k=p$, the Leibniz formula implies that $D^p$ is a derivation of $A$ so that the Lie algebra $\Der(A)$ is closed under the mapping $D\mapsto D^p$.  Investigating the relations between the operations of raising to the $p^{\rm th}$ power and the commutator bracket in the Lie algebra $\frak {gl}(A)=\End(A)$ for an associative algebra $A$ of positive characteristic leads one to the notion of a restricted Lie algebra.  These algebras, along with the corresponding representation theory, were first systematically studied by Jacobson in \cite{J}.  In particular, Jacobson defined the notion of the restricted universal enveloping algebra $U_{\rm res.}(\frak g)$ of a restricted Lie algebra $\frak g$, and showed that the category of restricted $\frak g$-modules is naturally isomorphic to the category of (unital) $U_{\rm res.}(\frak g)$-modules.  Unlike the case for ordinary Lie algebras, if $\frak g$ is a finite dimensional restricted Lie algebra, the algebra $U_{\rm res.}(\frak g)$ is also finite dimensional.  Many Lie algebras of positive characteristic admit a restricted Lie algebra structure. In \cite{Jb}, Jacobson shows that a necessary and sufficient condition for an ordinary Lie algebra of positive characteristic to admit the structure of a restricted Lie algebra is that all derivations of the form $(\ad g)^p$ be inner.  Such is the case for any Lie algebras with a non-degenerate Killing form for example.

Hochschild first considered the cohomology theory of restricted Lie algebras in \cite {H}.  Rather than constructing a cochain complex explicitly, he uses the canonical projection 
\[U(\frak g)\to U_{\rm res.}(\frak g),\]
where $U(\frak g)$ is the ordinary universal enveloping algebra, to induce a canonical map from restricted Lie algebra cohomology to ordinary Lie algebra cohomology for a given coefficient module $M$.  That is, this projection gives every restricted $\frak g$-module the structure of an ordinary $\frak g$-module so that a free resolution of the ground field by $U_{\rm res.}(\frak g)$ modules is a free resolution of $U(\frak g)$-modules.  Some investigations into the usual interpretations of low dimensional cohomology spaces are made using these mappings and exact sequence arguments.  This approach is insufficient in that the terms in the resolution used to define the restricted cohomology are too big to be effective for computational purposes.  One of the goals of the current dissertation is to improve on the results in \cite {H} by obtaining a smaller free resolution of the ground field by $U_{\rm res.}(\frak g)$-modules so that we obtain a cochain complex that is capable of making computations.  We achieve this goal up to dimension $p$ in the abelian case and partially in the non-abelian case.

The organization of this dissertation is as follows.  In chapter 2,  we briefly recall the theory of ordinary Lie algebra cohomology, including an explicit description of the Chevalley-Eilenberg cochain complex $\{C^*(\frak g;M),\delta\}$.  We then describe the general Cartan-Eilenberg theory of cohomology and show that the cohomology of the Chevalley-Eilenberg complex is isomorphic to the derived cohomology of a certain free resolution of the ground field by $U(\frak g)$-modules.  We then briefly recall some of the common algebraic interpretations of the low dimensional cohomology spaces of a Lie algebra including classes of extensions of modules and infinitesimal deformations of Lie algebras.  The second section of chapter 2 gives the definition of a restricted Lie algebra as well as reviews general theorems about the structure of these algebras and their representation theory.  In particular, we note that the restricted version of the Poincar\'e-Birkhoff-Witt theorem (Theorem (\ref {thm:ualg})) is valid so that there is a one-to-one correspondence between representations of a restricted Lie algebra and (unital) representations of its enveloping algebra $U_{\rm res.}(\frak g)$.  Chapter 3 contains our results on restricted Lie algebra cohomology.  In section (3.1), we construct an augmented complex 
\[C_*\stackrel{\epsilon}{\longrightarrow}\mathbb F\longrightarrow 0\]
of free $U_{\rm res.}(\frak g)$-modules, and show that this complex is exact in dimensions less than $p$ (Theorem (\ref {th:abelianresolution})). We then define the restricted cohomology of a restricted Lie algebra $\frak g$ with coefficients in the restricted module $M$ as
\[H^k({\frak g};M)=\Ext^k_{U_{\rm res.}({\scriptstyle \frak g})}(\mathbb F,M).\]
We also give the dimensions of the cochain spaces $C^k(\frak g;M)$ in the abelian case and show that the dimension of $C^k(\frak g)$ is the same as that of the symmetric algebra $S^k(\frak g)$ (Corollary (\ref {cor:dimensions})). In section (3.2), we explicitly construct the cochain spaces $C^k(\frak g;M)$ for $k\le 3$  and coboundary operators $\delta^k:C^k(\frak g;M)\to C^{k+1}(\frak g;M)$ for $k\le 2$ and show that the dimensions are the same as those in the abelian case.  Many proofs in this section are computational in nature, and the computations involve combinatorial identities modulo $p$ that are themselves interesting.  It is shown that there is no difference between ordinary cohomology and restricted cohomology in dimension 0, and Theorem (\ref {thm:1}) explicitly describes the restricted cohomology space $H^1(\frak g;M)$ as a subspace of the ordinary cohomology space $H^1_{\rm cl.}(\frak g;M)$.   We show that the canonical map
\[H^2(\frak g;M)\longrightarrow H^2_{\rm cl.}(\frak g;M)\]
is not injective in general, and we investigate its kernel.
In section (3.3), we develop the restricted analogs of the algebraic interpretations of low dimensional cohomology spaces and show that the cohomology defined by our complex encodes these notions.  In particular, we show that classes of restricted outer derivations of a restricted Lie algebra $\frak g$ coincide with $H^1(\frak g;\frak g)$ (Theorem (\ref {thm:outerderivations})); classes of restricted extensions of restricted modules $N$ by $M$ coincide with $H^1(\frak g;\Hom_{\mathbb F}(N,M))$ (Theorem (\ref {thm:rext})); classes of restricted extensions of restricted algebras $\frak g$ by $\frak h$ coincide with $H^2(\frak g;\frak h)$ (Theorem (\ref {thm:rextalg})); and classes of restricted infinitesimal deformations of restricted Lie algebras coincide with $H^2(\frak g;\frak g)$ (Theorem (\ref {thm:rinfinitesimaldeformations})).  Section (3.4) is a short discussion on the multiplicative structure of our complex. The dissertation concludes with a summary of the main results together with some remarks about further research.

\newpage
\pagestyle{myheadings}
\chapter{Background}
\thispagestyle{myheadings} 
\markright{  \rm \normalsize CHAPTER 2. \hspace{0.5cm}
  COHOMOLOGY OF RESTRICTED LIE ALGEBRAS }

\section{Lie Algebra Cohomology}

\subsection {The Chevalley-Eilenberg Cochain Complex} 
Historically, some of the first clues leading to the investigation of the cohomology groups of a finite dimensional Lie algebra were born out of attempts to generalize two lemmas belonging to Whitehead. The first is a result used in the proof of the complete reducibility of modules over a semi-simple Lie algebra. The second is the key result in the proof of a certain splitting theorem due to Levi.  To establish the flavor of Lie algebra cohomology theory, we give the precise statements here.

\begin{lemma}[Whitehead]
\label{whitehead1}
Suppose that ${\frak g}$ is a semi-simple Lie algebra over a field $\mathbb F$ of characteristic zero and suppose that $M$ is a ${\frak g}$-module. If $\phi:{\frak g}\to M$ is an $\mathbb F$-linear mapping satisfying 
\[ \phi([gh])=g\phi(h)-h\phi(g)\]
for all $g,h\in {\frak g}$, then there exists an element $m\in M$ such that $\phi(g)=gm$ for all $g\in {\frak g}$.
\qed
\end{lemma}

As we will see, this is nothing more than a statement about the triviality of the 1-dimensional cohomology groups of a semi-simple Lie algebra.  The second result concerns 2-dimensional cohomology groups.  The precise statement is as follows.

\begin{lemma}[Whitehead]
\label{whitehead2}
Suppose that ${\frak g}$ is a semi-simple Lie algebra over a field $\mathbb F$ of characteristic zero and suppose that $M$ is a finite dimensional ${\frak g}$-module. If $\phi:\Lambda^2{\frak g}\to M$ is a skew symmetric bilinear mapping satisfying 
\[ \phi([gh],f)+f\phi(h,g)+\phi([hf],g)+g\phi(f,h)+\phi([fg],h)+h\phi(g,f)=0,\]
then there exists a linear mapping $\psi:{\frak g}\to M$ such that 
\[ \phi(g,h)=g\psi(h)-h\psi(g)-\psi([gh]).\] 
\qed
\end{lemma}

Other sources for this theory include the study of the topology of Lie groups and vector fields on Lie groups by E. Cartan.  Lie algebras of smooth vector fields on manifolds are infinite dimensional and the cohomology theory described below must be suitably modified.  The details were completed  by Gelfand and Fuchs in \cite{FG1} and \cite {FG2}.  In this dissertation, we will mainly content ourselves with finite dimensional Lie algebras, although we put no restrictions on the characteristic of the underlying field $\mathbb F$  unless explicitly stated.  We begin with an explicit description of the standard complex used in the computation of Lie algebra cohomology.  It was originally constructed by Chevalley and Eilenberg in \cite {ChE}.  If ${\frak g}$ is a Lie algebra over $\mathbb F$ and $M$ is a ${\frak g}$-module, a  $q$-dimensional cochain of ${\frak g}$ with coefficients in $M$ is a skew-symmetric, $q$-linear map on ${\frak g}$ taking values in $M$.  The totality of all such maps comprises a vector space 
\[ C^q({\frak g};M)=\Hom_{\mathbb F}(\Lambda^q{\frak g},M) \]
over $\mathbb F$ under pointwise addition and scalar multiplication.  We set $C^q({\frak g};M)=0$ if $q<0$ and if $q=0$, we identify $C^0({\frak g};M)$ with $M\cong \Hom_{\mathbb F}(\mathbb F,M)$.  If $\phi\in C^q({\frak g};M)$, then $\phi$ determines an element $\delta\phi\in  C^{q+1}({\frak g};M)$ by the formula 
\begin{eqnarray*}
\delta\phi(g_1,\dots,g_{q+1})&=&\sum_{1\le s<t\le q+1}(-1)^{s+t-1}\phi([g_s,g_t],g_1,\dots,\widehat {g_s},\dots,\widehat {g_t},\dots,g_{q+1}) \\
 &+&\sum_{1\le s\le q+1}(-1)^s g_s\phi(g_1,\dots,\widehat {g_s},\dots,g_{q+1})
\end{eqnarray*}
where the symbol $\widehat {g_s}$ indicates that this term is to be omitted. It is easy to see that the mapping $\phi\mapsto \delta\phi$ is a linear transformation $\delta:C^q({\frak g};M)\to C^{q+1}({\frak g};M)$ and a direct verification shows that $\delta^2=0$.  Therefore $\{C^*({\frak g};M),\delta\}$ is a complex and its $q$th cohomology group is called the $q$-dimensional cohomology group (space) of ${\frak g}$ with coefficients in $M$ and is denoted by $H^q({\frak g};M)$.  We will denote the $q$-dimensional cocycles and coboundaries by $Z^q({\frak g};M)$ and $B^q({\frak g};M)$ respectively.  In dimensions $q=0,1$ and 2  the coboundary operator formula reduces to 
\begin{eqnarray*}
\delta(m)(g)&=&-gm\\
\delta\phi(g,h)&=&-g\phi(h)+h\phi(g)+\phi([gh])\\
\delta\phi(g,h,f)&=&\phi([gh],f)-\phi([gf],h)+\phi([hf],g) \\ 
 &-&g\phi(h,f)+h\phi(g,f)-f\phi(g,h)
\end{eqnarray*}
An examination of these formulae establishes the link between Whitehead's lemmas and Lie algebra cohomology. That is, both Lemmas \ref{whitehead1} and \ref{whitehead2} are consequences of the following theorem.
\newtheorem{vanish}[lemma]{Theorem}
\begin{vanish} 
\label{thm:vanish}
If ${\frak g}$ is a finite dimensional semi-simple Lie algebra over a characteristic zero field $\mathbb F$, then $H^1({\frak g};M)=0$ and $H^2({\frak g};M)=0$ for all finite dimensional ${\frak g}$-modules $M$.
\qed
\end{vanish}

We remark here that if $M=M_1\oplus M_2$, then easily we have 
\[H^q({\frak g};M)=H^q({\frak g};M_1)\oplus H^q({\frak g};M_2)\]
so that this result,  together with the complete reducibility of modules over semi-simple algebras of characteristic zero (Weyl's Theorem), reduces the computation of $H^q({\frak g};M)$ to the case $M$ irreducible.  Note that if ${\frak g}M=0$, we must have $\dim M=1$ by irreducibility so that $M$ is isomorphic to $\mathbb F$ and a $q$-cochain is a skew-symmetric $q$-linear form on ${\frak g}$ with values in $\mathbb F$.  Since ${\frak g}M=0$, the second term in the coboundary formula vanishes and we have
\[\delta\phi(g_1,\dots,g_{q+1})=\sum_{1\le s<t\le q+1}(-1)^{s+t-1}\phi([g_s,g_t],g_1,\dots,\widehat {g_s},\dots,\widehat {g_t},\dots,g_{q+1}).\] 
In the case of trivial coefficients, we usually shorten our notation and write $C^q({\frak g})$ in place of $C^q({\frak g};\mathbb F)$.  For semi-simple Lie algebras, the cohomology spaces with coefficients in $\mathbb F$ are the most interesting since they correspond to cohomology groups for Lie groups.  Indeed, if ${\frak g}$ is the Lie algebra of a Lie group $G$, then a $q$-cochain in $C^q({\frak g})$ gives rise to a right invariant differential form on $G$ so that we have inclusion $C^q({\frak g})\to \Omega^q(G)$ where $\Omega^*(G)$ denotes the de Rahm complex of the group $G$.  The trivial action ensures that this inclusion commutes with the differential so that we have inclusion of the complex $C^*({\frak g})$ into the de Rahm complex $\Omega^*(G)$.  If $G$ is compact, then this map induces an isomorphism in cohomology. The case ${\frak g}M\ne 0$ for semi-simple ${\frak g}$ and finite dimensional irreducible $M$ is not interesting because of the following theorem also due to Whitehead.
\begin{vanish}[Whitehead]
\label{thm:last}
If ${\frak g}$ is a finite dimensional semi-simple Lie algebra over a characteristic zero field $\mathbb F$, and $M$ is a finite dimensional irreducible ${\frak g}$-module such that ${\frak g}M\ne 0$, then $H^q({\frak g};M)=0$ for all $q\ge 0$.
\qed
\end{vanish}

Efforts to give algebraic proofs of these results of Whitehead as well as Weyl's theorem together provided the first clues for Chevalley and Eilenberg to give the preceding definition of $H^q({\frak g};M)$.


\subsection{Cartan-Eilenberg Definition of Cohomology}
In this subsection, we will give the Cartan-Eilenberg definition of the cohomology groups of a Lie algebra and show that the explicit definition given above computes this cohomology by means of a certain standard free resolution of $\mathbb F$.  The key notion is exploiting the correspondence between representations of a Lie algebra ${\frak g}$ and (unitary) representations of its universal enveloping algebra $U({\frak g})$.  Therefore we begin by recalling some facts concerning the universal enveloping algebra of a Lie algebra.  If ${\frak g}$ is a Lie algebra, we denote the tensor algebra by $T({\frak g})$.   Recall that as a vector space,  $T({\frak g})$ is given by 
\[T({\frak g})=\bigoplus_{n=0}^\infty (\underbrace {{\frak g}\otimes\cdots\otimes{\frak g}}_n) \]
and the multiplication is defined on homogeneous generators by juxtaposition.
We denote by $I$ the two-sided ideal in $T({\frak g})$ generated by all elements of the form
\[g\otimes h-h\otimes g-[gh]\]
and define the universal enveloping algebra for ${\frak g}$ as the quotient $U({\frak g})=T({\frak g})/I$.  We write the image of a generator $g_1\otimes\cdots\otimes g_n\in T({\frak g})$ as $g_1\cdots g_n\in U({\frak g})$.  The canonical augmentation $\epsilon:T({\frak g})\to \mathbb F$ vanishes on $I$ and hence we have an augmentation of the algebra $U({\frak g})$ which we also denote by $\epsilon$.  We denote the kernel of $\epsilon:U({\frak g})\to \mathbb F$ by $U({\frak g})^+$.  The importance of $U({\frak g})$ in our cohomology theory lies in the fact that there is a one-to-one correspondence between Lie algebra representations of ${\frak g}$ (${\frak g}$-modules) and unitary representations of $U({\frak g})$.  Consequently, the cohomology theory of a Lie algebra can be entirely constructed using the associative algebra $U({\frak g})$. Indeed, if we regard $\mathbb F$ as a trivial $U({\frak g})$-module (that is $a\lambda=\epsilon(a)\cdot\lambda$ for all $a\in U({\frak g})$ and all $\lambda\in \mathbb F$), the cohomology spaces $H^q({\frak g};M)$ defined above coincide with the spaces $\Ext^q_{U({\scriptstyle\frak g})}(\mathbb F,M)$ where $M$ is regarded as a unitary $U({\frak g})$-module in the natural fashion. 
To see this, let us define $D^q$ as the space of $q$-linear (non-alternating) forms on $U({\frak g})^+$ with values in $M$ and a coboundary operator $\delta:D^q\to D^{q+1}$ by the formula
\[\delta f(x_1,\dots,x_{q+1})=x_1f(x_2,\dots,x_{q+1})+\sum_{j=2}^{q+1}(-1)^{j-1}f(x_1,\dots,x_{j-1}x_j,\dots,x_{q+1}).\]
We remark that this formula is a special case of the cohomology theory of associative algebras as defined by Hochschild in \cite {H2} if we set the appropriate operations to zero because of the trivial action on $\mathbb F$.  It is possible to give an explicit cochain map which establishes a natural isomorphism between the cohomology of the complex $\{D^*,\delta\}$ and the cohomology of the Chevalley-Eilenberg complex $\{C^*({\frak g};M),\delta\}$.  Briefly, for every $q$-cochain $f\in D^q$, we define a cochain $f^\prime\in C^q({\frak g};M)$ by the formula
\[f^\prime(g_1,\dots,g_n)=\sum_{\sigma\in S_q}\sgn\sigma f(g_{\sigma(1)},\dots,g_{\sigma(q)}).\]
A direct verification shows that $(\delta f)^\prime=\delta f^\prime$ so that the assignment $f\mapsto f^\prime$ induces a map of cohomology groups.  To see that this induced map is actually an isomorphism, one shows that both of the above complexes are the equivariant cohomology groups derived from free resolutions of $\mathbb F$ by $U({\frak g})$-modules.  The general theory of algebraic complexes then implies that the resulting cohomology groups are naturally isomorphic.    

In order to better motivate some of our later manipulations, we briefly recount the constructions of these resolutions here.  The cochain complex $\{C^*({\frak g};M),\delta\}$ is obtained from the following resolution.  We set $\C_0=U({\frak g})$ regarded as a regular (left) $U({\frak g})$-module and we use the augmentation $\epsilon:\C_0\to \mathbb F$ as defined above.  For $q>0$, we set
\[\C_q=U({\frak g})\otimes \Lambda^q{\frak g}\]
with the natural $U({\frak g})$-module structure.  The boundary operator $d:\C_q\to \C_{q-1}$ is defined by the formula
\begin{eqnarray*}
d(x\otimes g_1\wedge\cdots\wedge g_q)&=&\sum_{j=1}^n(-1)^{j-1} xg_j\otimes g_1\wedge\cdots \widehat {g_j}\cdots\wedge g_q\\
 &+&\sum_{1\le s<t\le q}(-1)^{s+t-1} x\otimes [g_s,g_t]\wedge x_1\wedge\cdots \widehat {g_s}\cdots\widehat {g_t}\cdots\wedge g_q.
\end{eqnarray*}
Clearly each $\C_q$ is a free $U({\frak g})$-module and a direct computation shows that $d^2=0$.  The verification of the acyclicity of this complex is not so straightforward.  One must introduce a grading of the complex $\C_*$ and use this grading to introduce an increasing filtration of sub-complexes whose successive quotients are easily described in terms of the symmetric and alternating algebras on ${\frak g}$. These quotient complexes can easily be shown to be acyclic and it then follows from general homological algebra theorems that the entire complex $\C_*$ is itself acyclic. A detailed proof by Kozul using spectral sequences can be found in  \cite {W}, Theorem (7.7.2).  Now if $M$ is a unitary $U({\frak g})$-module, it is clear that 
\[C^q({\frak g};M)=\Hom_{U({\scriptstyle \frak g})}(\C_q, M)\] 
and the coboundary operator $\delta:C^q({\frak g};M)\to C^{q+1}({\frak g};M)$ is the dual of the boundary map $d$.  Therefore our original description of the cohomology of a Lie algebra is the equivariant cohomology group derived from $M$ and the free resolution $\C_*\to \mathbb F$.
The construction in the associative case is less complicated.  Here we set $\D_0=U({\frak g})$ and use the same augmentation. For $q>0$ we define
\[ \D_q=U({\frak g})\otimes\underbrace{U({\frak g})^+\otimes\cdots\otimes U({\frak g})^+}_q \]
and the boundary operator is given by
\begin{eqnarray*}
d(x\otimes g_1\otimes\cdots\otimes g_q)&=&xg_1\otimes g_2\otimes\cdots\otimes g_q\\
 &+&\sum_{j=2}^n(-1)^{j-1} x\otimes g_1\otimes\cdots\otimes g_{j-1}g_j\otimes\cdots \otimes g_q.
\end{eqnarray*}
In this case we can define an explicit chain homotopy operator $H:\D_q\to \D_{q+1}$ that shows the complex $\D_*$ is acyclic.  Namely, we define $H:\mathbb F\to \D_0$ to be the inclusion of $\mathbb F$ into $U({\frak g})$, and on $\D_1$, we define $H$ by the formula $H:x\mapsto 1\otimes (x-\epsilon(x))$.  For $q>0$, we define $H:\D_q\to \D_{q+1}$ by
\[ H(x\otimes g_1\otimes \cdots\otimes g_q)=1\otimes(x-\epsilon(x))\otimes g_1\otimes \cdots\otimes g_q.\]
It is easy to see that $dH+Hd$ is the identity map on the complex $\D_*$ so that $\D_*$ is acyclic.  Moreover, we evidently have
\[D^q=\Hom_{U({\scriptstyle \frak g})}(\D_q, M)\]
and the coboundary operator $\delta:D^q\to D^{q+1}$ is the dual of the boundary map $d$. Therefore the cohomology of the associative complex $\{D^*,d\}$ is the equivariant cohomology group derived from $M$ and the free resolution $\D_*\to \mathbb F$.  

To complete the argument that our two cohomology groups are isomorphic, we recall that if $A$ is an associative algebra and $\C_*\to F$ and $\D_*\to F$ are two resolutions of $A$-modules of the same $A$-module $F$, and if $\C_*$ is free over $A$ while $\D_*$ is acyclic, then the identity map $F\to F$ can be extended to a chain map $\C_*\to \D_*$.  Moreover, any two such extensions are chain homotopic.  It follows that there is a canonical induced homomorphism $H(\C_*)\to H(\D_*)$.  In our case, we can interchange the roles of $\C_*$ and $\D_*$ to see that this induced homomorphism is actually an isomorphism.  The argument is completed by noting that our map $f\mapsto f^\prime$ of the complex $D^*$ into $C^*$ is the dual of a certain chain map $\C_*\to \D_*$ and hence induces an isomorphism in cohomology.  

We conclude this subsection with a final remark on the acyclicity of the complex $\C_*$ which may be useful in the sequel.  Quite generally, if $A$ is an associative unital algebra over $\mathbb F$, $B\subset A$ is a subalgebra of $A$, and ${}_BM$ is a left $B$-module, then the abelian group $A_B\otimes_B {}_BM$ is canonically a left $A$-module with the action
\[\alpha(a\otimes m)=(\alpha a)\otimes m.\]
where $\alpha,a\in A$ and $ m\in M$. We call this $A$-module the induced module and write $\induced_B^A(M)$.  If we define the homology of $B$ with coefficients in $M$ as $H_q(B;M)=\Tor^{B}_q(\mathbb F;M)$, then the acyclicity of the complex $\C_*$ follows immediately from the following result known as Shapiro's Lemma.
\begin{lemma}[Shapiro]
\label{thm:shapiro}
$H_q(B;M)$ is isomorphic to $H_q(A;\induced_B^A(M))$ for all $q\ge 0$.
\end{lemma}
If we apply the above considerations to $U({\frak g})$ and $U({\frak h})$, where ${\frak h}\subset {\frak g}$ is a subalgebra of ${\frak g}$, we get an $U({\frak g})$-module $\induced_{U({\frak h})}^{U({\frak g})}(M)=U({\frak g})\otimes_{U({\scriptstyle\frak h})} M$.  In particular, if ${\frak h}=0$ so that a ${\frak h}$-module $M$ is just a vector space; we can take $M=\mathbb F$ and we have $\induced_0^{U({\frak g})}(\mathbb F)=U({\frak g})\otimes_{\mathbb F} \mathbb F=U({\frak g})$.  Therefore Shapiro's lemma states
\[H_q({\frak g};U({\frak g}))\cong H_q(\mathbb F;\mathbb F)\]
and the later spaces are easily seen to vanish for all $q$. The argument is complete upon noticing that 
\[C_q({\frak g};U({\frak g}))=\C_q\otimes_{U({\frak g})}U({\frak g})=\C_q.\]

Having this definition of the cohomology spaces, we now turn our attention toward methods of computing these spaces for a given Lie algebra and coefficient module.


\subsection{Hochschild-Serre Spectral Sequence}
In this subsection we describe one of the main computational tools in Lie algebra cohomology - the Hochschild-Serre spectral sequence.  In 1953, Hochschild and Serre studied the algebraic relations in the cohomology of a group that arise from group extensions $N\to G\to G/N$ by giving a filtration of the complex of cochains of $G$ with coefficients in an certain $G$-module and studying the resulting spectral sequence.  In a follow up paper \cite {HS}, they gave the analogous filtration of the complex $\{C^*({\frak g};M),\delta\}$ and showed that the corresponding spectral sequence abuts to the cohomology of the Lie algebra.  For the convenience of the reader, we briefly recall the details of the construction of the spectral sequence of a filtered complex here.  We recall that a (degree $+1$) complex (of abelian groups, $A$-modules, etc.) $K=\{K^q,d\}$ is filtered if we are given a decreasing sequence of subcomplexes
\[K=F^0K\supset F^1K\supset\cdots\supset F^nK\supset F^{n+1}K=0.\]
Equivalently, for each $q\ge 0$, we have a decreasing filtration
\[K^q=F^0K^q\supset F^1K^q\supset\cdots\supset F^{n_q}K^q\supset F^{n_q+1}K^q=0\]
of $K^q$ such that $d(F^sK^q)\subset F^sK^{q+1}$ for all $s$ and all $q$. To each filtered complex $K$ we have an associated graded complex 
\[\Gr K=\bigoplus_{s\ge 0} \Gr^s K\]
where $\Gr^s K=F^sK/F^{s+1}K$ and the differential $d:\Gr^s K\to \Gr^{s+1} K$ is induced by $d:K\to K$.  The inclusion $F^sK\to K$ induces a map $H(F^sK)\to H(K)$ and we let 
\[F^sH(K)=\Im(H(F^sK)\to H(K)).\]  
We therefore have a decreasing filtration
\[H(K)=F^0H(K)\supset F^1H(K)\supset\cdots\supset F^nH(K)\supset F^{n+1}H(K)=0\]
and an associated graded complex
 \[\Gr H(K)=\bigoplus_{s\ge 0} \Gr^s H(K).\]
 If we keep track of the grading in $K$, we have 
 \[\Gr H(K)=\bigoplus_{s,q\ge 0} \Gr^s H^q(K).\]
The Hochschild-Serre spectral sequence is a special case of the following general theorem.
\begin{vanish}
\label{thm:specseq}
If $K$ be a filtered complex, then there exists a spectral sequence $\{E^{s,q}_r\}$, $s,q,r\ge 0$, with 
\begin{eqnarray*}
E^{s,q}_0&=&F^sK^{s+q}/F^{s+1}K^{s+q}  \\
E^{s,q}_1&=&H^{s+q}(\Gr^s K) \\
E^{s,q}_\infty&=&\Gr^s(H^{s+q}(K)).
\end{eqnarray*}
\end{vanish}
The last relation is written $E_r\Rightarrow H(K)$ and the spectral sequence is said to abut to $H(K)$.  

Now, if ${\frak g}$ is a Lie algebra, ${\frak h}\subset {\frak g}$ is a subalgebra and $M$ is a ${\frak g}$-module,  Hochschild and Serre defined a filtration of the standard complex $\{C^*({\frak g};M),\delta\}$ in \cite {HS} as follows.  Define $F^sC^q=F^sC^q({\frak g};M)\subset C^q({\frak g};M)$ by 
\[F^sC^q=\{\phi\in C^q({\frak g};M) : \phi(g_1,\dots,g_q)=0\ \hbox {\rm whenever $q-s+1$ of the $g_i$ lie in ${\frak h}$}\}.\]
One can easily see that $\delta(F^sC^q)\subset F^sC^{q+1}$ so that we have the filtration
\[C^q({\frak g};M)=F^0C^q\supset F^1C^q\supset F^2C^q\supset\cdots\supset F^qC^q\supset F^{q+1}C^q=0\]
which is compatible with the differential $\delta$.  We combine the main results in \cite {HS} in the following theorem.
\begin{vanish}
\label{thm:hochserre}
If ${\frak g}$ is a Lie algebra and ${\frak h}$ is a subalgebra of ${\frak g}$, then there exists a spectral sequence $\{E^{s,q}_r\}$ such that:
\begin{itemize}
\item[(\romannumeral 1)] $E^{s,q}_0=C^q({\frak h};\Hom_{\mathbb F}(\Lambda^s ({\frak g/h});M))$.
\item[(\romannumeral 2)] The differential $d^{s,q}_0:E^{s,q}_0\to E^{s,q+1}_0$ is the usual differential 
\[ \delta:C^q({\frak h};\Hom_{\mathbb F}(\Lambda^s ({\frak g/h});M))\to C^{q+1}({\frak h};\Hom_{\mathbb F}(\Lambda^s ({\frak g/h});M)), \]
so that $E^{s,q}_1=H^q({\frak h};\Hom_{\mathbb F}(\Lambda^s ({\frak g/h});M))$.
\item[(\romannumeral 3)] If ${\frak h}$ is an ideal, then $E^{s,q}_2=H^s({\frak g/h};H^q({\frak h};M))$.
\item[(\romannumeral 4)] $E^{s,q}_r\Rightarrow H^{s+q}({\frak g};M)$.
\end{itemize}
\end{vanish}


\subsection{Algebraic Interpretations}
We close this section with some remarks on general algebraic interpretations of low dimensional cohomology groups of a Lie algebra. The corresponding notions will serve as motivation for the sequel.

First, for any coefficient module $M$, the space $H^0({\frak g};M)$ is naturally isomorphic to the space $M^{\frak g}$ of ${\frak g}$-invariants. Recall that an element $m\in M$ is a ${\frak g}$-invariant if $gm=0$ for all $g\in {\frak g}$.  This isomorphism is easy to understand given that $H^0({\frak g};M)=Z^0({\frak g};M)$ and $\delta (m)=0$ if and only if $-gm=0$ for all $g$.  

For a slightly more interesting example, recall that a derivation of a Lie algebra $\frak g$ is a linear map $D:\frak g\to \frak g$ satisfying
\[D[gh]=[gD(h)]+[D(g)h]\]
for all $g,h\in \frak g$. The Lie commutator of two derivations is again a derivation so that the space $\Der(\frak g)$ of all derivations of $\frak g$ is a Lie subalgebra of $\frak {gl}(\frak g)$.  If $g\in\frak g$ is fixed, then the Jacobi identity implies that the map $\ad g:\frak g\to \frak g$ defined by $\ad g(h)=[gh]$ is a derivation of $\frak g$.  Such a derivation is called inner.  In fact,  the map $\ad :\frak g\to \Der (\frak g)$ is a Lie algebra homomorphism and $\ad(\frak g)$ is an ideal in $\Der (\frak g)$.  By definition, an outer derivation of $\frak g$ is an element of the quotient $\Der (\frak g)/\ad (\frak g)$.  Now, the map $\ad :\frak g\to \Der (\frak g)$ gives $\frak g$ the structure of a $\frak g$-module and the 1-dimensional cochain space with coefficients in $\frak g$ is $C^1(\frak g;\frak g)=\Hom_{\mathbb F}(\frak g,\frak g)$.  If $\phi\in C^1(\frak g;\frak g)$, then 
\begin{eqnarray*}
\delta\phi(g,h)=0&\Longleftrightarrow& \phi([gh])-[g\phi(h)]+[h\phi(g)]=0 \\
 &\Longleftrightarrow& \phi([gh])=[g\phi(h)]+[\phi(g)h].
\end{eqnarray*}
Therefore $\phi$ is a cocycle if and only if $\phi$ is a derivation of $\frak g$.  Moreover, we have $C^0(\frak g;\frak g)=\frak g$ and 
\[\delta g(h)=-[hg]=[gh]=\ad g (h)\]
so that $\Im\delta=\ad(\frak g)$.  It follows that $H^1(\frak g;\frak g)=\Der(\frak g)/\ad (\frak g)$.  That is $H^1(\frak g;\frak g)$ is canonically isomorphic to the space of outer derivations on $\frak g$.

The space $H^1(\frak g;\frak g)$ has another interpretation. Recall that a  one dimensional right extension of $\frak g$ is by definition a short exact sequence of Lie algebras and their homomorphisms
\begin{equation}
\label{rightsequence}
0\to\frak g \to \frak g'\to \mathbb F\to 0
\end{equation}
with the Lie bracket in $\frak g'=\frak g\oplus \mathbb F$ defined by
\begin{eqnarray}
\label{rightbracket}
[(g_1,\lambda_1),(g_2,\lambda_2)]&=&([g_1,g_2]-\lambda_1c(g_2)+\lambda_2c(g_1),0)
\end{eqnarray}
where $c:\frak g\to \frak g$ is a linear map. Two one dimensional right extensions of $L$ are equivalent if they can be included in a commutative diagram
\begin{equation} 
\label{diagram}
\begin{array}{ccccccccc}
 0 & \longrightarrow & \frak g & \longrightarrow& \frak g'& \longrightarrow &\mathbb F & \longrightarrow &0 \\
&&\|&\circ&\downarrow&\circ&\|& \\
0&\longrightarrow&\frak g&\longrightarrow&\frak g''&\longrightarrow&\mathbb F&\longrightarrow&0 
\end{array}
\end{equation}
It can be shown that $c\in C^1(\frak g;\frak g)$ is a cocycle if and only if  the bracket in (\ref {rightbracket}) satisfies the Jacobi identity (it is clearly bilinear and skew symmetric).  If you write out the left hand side  of the Jacobi identity for (\ref {rightbracket}) using $(g,\alpha),(h,\beta),(f,\gamma)\in \frak g\oplus \mathbb F$, it simplifies to:
\[([g[hf]]+[h[fg]]+[f[gh]]+\alpha \delta c(h,f)+\beta \delta c(f,g)+\gamma \delta c(g,h),0)\]
where $\delta c$ is the coboundary of the cochain $c$.
Therefore if $c$ is a cocycle, then (\ref {rightbracket}) defines a Lie bracket on $\frak g'$.  If we assign to each cocycle $c\in C^1(\frak g;\frak g)$ the sequence (\ref {rightsequence}) with bracket (\ref {rightbracket}), then our previous remark shows that this assignment is surjective.  If $c'=c+\delta b$ for some $b\in C^0(\frak g;\frak g)=\frak g$, then the corresponding sequences may be included in a commutative diagram like that above whose middle vertical map is defined by $(g,\lambda)\mapsto (g+\lambda b,\lambda)$.  This map is a Lie algebra homomorphism making the diagram commute so that our assignment is well defined.  Finally, one can show that equivalent sequences come from cohomologous cocycles so that the space $H^1(\frak g;\frak g)$ can be naturally identified with equivalence classes of one dimensional right extensions of $\frak g$.  More generally, if $M$ is any $\frak g$-module, then $H^1(\frak g;M)$ is naturally isomorphic to the space of equivalence classes of 1-dimensional right extensions of $M$.  By definition, such an extension is an exact sequence of $\frak g$-modules
\begin{equation}
0\to M\to M'\to \mathbb F\to 0
\end{equation}
where $\mathbb F$ is considered as a trivial $\frak g$-module.  Choosing a preimage for $1\in\mathbb F$ determines a linear map $\frak g\to M$.  The module condition on $M$ makes this map a cocycle and different choices of preimages give cohomologous cocycles.

For another example, recall that a central extension of a Lie algebra $\frak g$ is an exact sequence of Lie algebras
\begin{equation}
\label{centralsequence}
0\to \mathbb F\to \frak g'\to \frak g\to 0
\end{equation}
with the Lie bracket in $\frak g'=\mathbb F\oplus \frak g$ given by 
\begin{equation}
\label{centralbracket}
[(\lambda_1,x_1),(\lambda_2,x_2)]=(c(x_1,x_2),[x_1,x_2])
\end{equation}
where $c:\frak g\times\frak g\to \mathbb F$ is a skew-symmetric bilinear map.  Two central extensions are equivalent if they can be included in a diagram analogous to (\ref {diagram}).  Here, the Jacobi identity for the bracket in (\ref {centralbracket}) is equivalent to $c\in C^2(\frak g)$ being a cocycle.  If $c'=c+\delta b$ for some $b\in C^1(\frak g)=\frak g^*$, then the map $(\lambda,x)\mapsto (\lambda+b(x),x)$ is an equivalence of the corresponding central extensions (\ref {centralsequence}) so that $H^2(\frak g)$ is naturally isomorphic to the space of equivalence classes of central extensions of $\frak g$.  

We close this section with a cohomological description of infinitesimal deformations of a Lie algebra $\frak g$.  Recall that an infinitesimal deformation of $\frak g$ is a map $\eta:\frak g\times\frak g\times \mathbb F\to \frak g$ written
\begin{equation}
\label{deformation}
\eta(g,h,t)=[gh]_t=[gh]+tc(g,h)
\end{equation}
where $c\in C^2(\frak g,\frak g)$ and $[gh]_t$ is a Lie bracket $\pmod {t^2}$ in $\frak g$ for all $t\in \mathbb F$. If we write out the Jacobi identity for the bracket (\ref {deformation}) and simplify, we have 
\begin{equation}
\label{defjacobi}
t\left (c(g,[hf])+[gc(h,f)]+c(h,[fg])+[hc(f,g)]+c(f,[gh])+[fc(g,h)]\right ).
\end{equation} 
Of course we recognize the parenthesized term in (\ref {defjacobi}) as $-\delta c(g,h,f)$ so that $\eta$ is an infinitesimal deformation if and only if $c\in C^2(\frak g,\frak g)$ is a 2-cocycle.  Two infinitesimal deformations $\eta_1$ and $\eta_2$ are equivalent if there is a linear map $\xi:\frak g\to \frak g$ such that 
\begin{equation}
\label{defequivalence}
c_1(g,h)=c_2(g,h)+[g\xi(h)]+[\xi(g)h]-\xi([gh])
\end{equation}
for all $g,h\in \frak g$.  With this definition, we see from (\ref {defequivalence}) that two infinitesimal deformations $\eta_1$ and $\eta_2$ are equivalent if and only if the corresponding cocycles $c_1$ and $c_2$ are cohomologous.  Therefore we can identify the space of equivalence classes of infinitesimal deformations with $H^2(\frak g,\frak g)$.


\section {Restricted Lie Algebras}
\subsection {Guiding Examples and Definitions}
Lie algebras over fields of characteristic $p>0$ often posses an additional structure involving a set map ${\frak g}\to {\frak g}$. These objects were first systematically studied by Jacobson in \cite {J} where he termed them restricted Lie algebras.  The corresponding cohomology theory was first examined by Hochschild in \cite {H}.  In order to better motivate the formal definitions, we begin by describing two examples that serve as prototypes for the general notion of a restricted Lie algebra. First suppose that $A$ is an associative algebra over $\mathbb F$ and recall that a derivation of $A$ is a linear map $D:A\to A$ that satisfies the Leibniz rule
\[D(ab)=aD(b)+D(a)b\]
for all $a,b\in A$.  The Lie commutator $[D,D^\prime]=DD^\prime-D^\prime D$ of two derivations is again a derivation so that the subspace $\Der A$ of linear transformations on $A$ that consists of all derivations on $A$ is a Lie subalgebra of the Lie algebra $\frak {gl}(A)$ of all linear transformations on $A$.  The composition of two derivations is, in general, not a derivation.  However, an easy induction argument shows that we do have the following more general version of the Leibniz rule.  Specifically, if $D\in \Der A$, $a,b\in A$ and $k\ge 1$, we have
\begin{eqnarray}
\label{leibniz}
D^k(ab)&=&\sum_{j=0}^k \bino {k}{j}D^j(a)D^{k-j}(b).
\end{eqnarray}
If we assume that the characteristic of the ground field $\mathbb F$ is positive, and we take $k=p=\mathop{\rm char}\mathbb F$, then $\bino {k}{j}=0\pmod {p}$ unless $j=0$ or $j=p$ so that equation (\ref{leibniz}) reduces to
\[D^p(ab)=aD^p(b)+D^p(a)b.\]
Therefore if $\mathop{\rm char}\mathbb F=p>0$, the Lie algebra $\Der A$ is closed under an additional mapping $D\mapsto D^p$.  For the second example, recall that every associative algebra $A$ has an underlying Lie algebra structure with Lie bracket $[a,b]=ab-ba$. We denote this Lie algebra by the corresponding Gothic letter ${\frak a}$.  If $A$ has an anti-automorphism $a\mapsto \overline a$, then the subset ${\frak g}$ of ${\frak a}$ that consists of all skew elements with respect to this anti-automorphism (i.e. those $a\in A$ satisfying $\overline a=-a$) is a Lie subalgebra of ${\frak a}$.  Moreover, if $\mathop{\rm char}\mathbb F=p>0$, then $\overline {a^p}=\overline a^p=(-a)^p=-a^p$ so that ${\frak g}$ is also closed under the mapping $a\mapsto a^p$. Investigating the algebraic relations between the operations of raising to the $p{\rm th}$ power and the Lie bracket in the underlying Lie algebra of an associative algebra leads to the abstract definition of a restricted Lie algebra. Before giving the definition, we recall that if $\mathbb F$ is a field of characteristic $p>0$ and $X$ and $Y$ are two (commuting) indeterminants, then we have the following well known relations in the polynomial ring $\mathbb F [X,Y]$:
\begin{eqnarray*}
(X-Y)^p&=&X^p-Y^p \\ 
(X-Y)^{p-1}&=& \sum_{j=0}^{p-1} X^jY^{p-1-j}.
\end{eqnarray*}
These relations imply corresponding relations for any two commuting elements $x$ and $y$ in an associative $\mathbb F$-algebra $A$.  In particular, if $b\in A$ is fixed and we denote the left and right multiplications by $b$ by $\lambda_b$ and $\rho_b$ respectively,  then we have
\begin{eqnarray*}
(\rho_b-\lambda_b)^p&=&\rho_b^p-\lambda_b^p=(\rho_b)^p-(\lambda_b)^p \\ 
(\rho_b-\lambda_b)^{p-1}&=& \sum_{j=0}^{p-1} \rho_b^j\lambda_b^{p-1-j}=\sum_{j=0}^{p-1} (\rho_b)^j(\lambda_b)^{p-1-j},
\end{eqnarray*}
or equivalently, for all $a,b\in A$, 
\begin{eqnarray}
\label{eq:relations}
\mbox{}[\cdots[[a\overbrace{b]b]\cdots b}^{p}]&=&[ab^p] \nonumber \\
\mbox{}[\cdots[[a\overbrace{b]b]\cdots b}^{p-1}]&=&\sum_{j=0}^{p-1} b^{p-1-j}ab^j
\end{eqnarray}
It is clear that $(\alpha a)^p=\alpha^p a^p$ for all $\alpha\in\mathbb F$ and $a\in A$.  Moreover, one can use the relation (\ref {eq:relations}) to show that 
\[(a+b)^p=a^p+b^p+\sum_{j=1}^{p-1} s_j(a,b)\]
where $js_j(a,b)$ is the coefficient of $X^{j-1}$ in $(\ad (Xa+b))^{p-1}(a)$, $X$ an indeterminate.  All of these considerations lead to the following definition.  
\begin{definition}
\label{def:maindef}
A restricted Lie algebra of characteristic $p>0$ is a Lie algebra ${\frak g}$ of characteristic $p$ together with a map ${\frak g}\to {\frak g}$, denoted by $g\mapsto g^{[p]}$, that satisfies
\begin{itemize}
\item[R1] $(\lambda g)^{[p]}=\lambda^pg^{[p]}$, 
\item[R2] $\ds (g+h)^{[p]}=g^{[p]}+h^{[p]}+ \sum_{\stackrel{g_j=g\ {\rm or}\  h}{\scriptscriptstyle g_1=g,g_2=h}} \frac{1}{\#(g)} [[[\cdots[[g_1g_2]g_3]\cdots]g_{p-1}]g_p]$,
\item[] where $\#(g)$ denotes the number of $g$'s among the $g_j$,
\item[R3] $[gh^{[p]}]=[[\cdots [g\underbrace{h]h]\cdots]h}_{p}].$
\end{itemize}
for all $g,h\in {\frak g}$ and all $\lambda\in\mathbb F$. 
\end{definition}
The map $g\mapsto g^{[p]}$ is referred to as the $p$-operator.  For notational ease, we will write all multiple Lie brackets with the notation
\[ [[[\cdots[[g_1g_2]g_3]\cdots]g_{p-1}]g_k]=[g_1,g_2,\dots,g_k] \]
so that in this notation we have
\begin{itemize} 
\item[{\it R2}] $\ds (g+h)^{[p]}=g^{[p]}+h^{[p]}+ \sum_{\stackrel{g_j=g\ {\rm or}\  h}{\scriptscriptstyle g_1=g,g_2=h}} \frac{1}{\#(g)}[g_1,g_2,\dots,g_p] $,
\item[{\it R3}] $[gh^{[p]}]=[g,\underbrace{h,\dots,h}_{p}].$
\end{itemize}
The remarks at the beginning of this subsection imply that if $A$ is an associative algebra over $\mathbb F$, then the underlying Lie algebra ${\frak a}$ of $A$ is a restricted Lie algebra with the operation $a^{[p]}=a^p$.  In particular, if $M$ is a vector space over $\mathbb F$, then the algebra $\End_{\mathbb F} M$ of $\mathbb F$-linear transformations $M\to M$ is a restricted Lie algebra which will be denoted ${\frak {gl}}(M)$.  A homomorphism $\phi$ from one restricted Lie algebra to another is a Lie algebra homomorphism $\phi:{\frak g}\to {\frak h}$ such that $\phi(g^{[p]})=\phi(g)^{[p]}$ for all $g\in {\frak g}$.  Restricted subalgebras, kernels and ideals are all defined in the obvious way.  In \cite {Jb}, Jacobson gives a necessary and sufficient condition in which an ordinary Lie algebra of characteristic $p>0$ admits the structure of a restricted Lie algebra.  Indeed, condition (R3) in Definition (\ref {def:maindef}) makes it clear that a necessary condition is that for every $h\in {\frak g}$, the derivation $(\ad h)^p$ is inner.  In \cite {Jb}, it is shown that this condition is also sufficient.  In fact, it suffices that $(\ad e_j)^p$ is inner for all $e_j$ in some basis for ${\frak g}$.  We remark that if ${\frak g}$ is restricted with respect to two $p$-operators $g\mapsto g^{[p]_1}$ and $g\mapsto g^{[p]_2}$, then the map 
\[f:g\mapsto g^{[p]_1}-g^{[p]_2}\] 
maps ${\frak g}$ into the center of ${\frak g}$.  Moreover, $f$ is $p$-semi-linear in the sense that $f(g+h)=f(g)+f(h)$ and $f(\lambda g)=\lambda^p f(g)$ for all $g,h\in {\frak g}$ and all $\lambda\in\mathbb F$.  The kernel of a $p$-semi-linear map is a subspace so that if two $p$-operators that make ${\frak g}$  a restricted Lie algebra agree on a basis, they are identical.  If $\frak g$ is finite dimensional with a non-degenerate Killing form, then every derivation is inner so that $\frak g$ admits a restricted Lie algebra structure.  Moreover, it is clear that in this case that the center of $\frak g$ is $0$ so that the $p$-operator on $\frak g$ is unique. For an example, let $\mathbb Z_p$ denote the cyclic group of order $p$ and let $A=\mathbb F[\mathbb Z_p]$ denote the group algebra of $\mathbb Z_p$ over $\mathbb F$.  As an $\mathbb F$-algebra, $A$ has a basis $\{1,x,x^2,\dots,x^{p-1}\}, (x^p=1)$ where $x\in\mathbb Z_p$ denotes a generator.  It can be shown that the derivation algebra $\Der A$ has a basis $D_j, j=0,1,\dots, p-1$, where 
\[D_j(x)=x^{j+1}.\]
Moreover, it is easy to verify that 
\begin{eqnarray*}
[D_i,D_j]&=&(i-j)D_{i+j} \\
D_0^p&=&D_0 \\
D_j^p&=&0\ (j>0).
\end{eqnarray*}
Finally, we remark that $\Der A$ is simple as a Lie algebra so that the $p$-operator $D^{[p]}=D^p$ is the only map $\Der A\to \Der A$ giving $\Der A$ a restricted Lie algebra structure. The algebra $\Der A$ is usually referred to as the Witt algebra.

We close this subsection by remarking that in the case when the Lie algebra $\frak g$ is abelian, that is $[\frak g\frak g]=0$, the $p$-operator is a map $\frak g\to \frak g$ satisfying $(\lambda g)^{[p]}=\lambda^pg^{[p]}$ and $(g+h)^{[p]}=g^{[p]}+h^{[p]}$ for all $g,h\in\frak g$ and $\lambda\in \mathbb F$.  A map with these properties is called $p$ semi-linear.  Recall that since $\mathbb F$ has characteristic $p>0$, the Frobenius map $\alpha: \lambda \mapsto \lambda^p$ is an automorphism of $\mathbb F$.   If $V$ is an abelian group, an $\mathbb F$-vector space structure on $V$ is completely determined by giving a ring homomorphism $\mathbb F\to \End(V)$ where $\End(V)$ denotes the ring of group endomorphisms $V\to V$.  Therefore if  $V$ is a vector space over $\mathbb F$, then the composition
\[\mathbb F\stackrel{\alpha^{-1}}{\longrightarrow}\mathbb F\longrightarrow\End(V)\]
gives another vector space structure on $V$ which we will denote by  $\overline V$. Now if $\mathbb F_p\subset\mathbb F$ denotes the cyclic subfield of $\mathbb F$ of order $p$, a result often referred to as Fermat's Little Theorem implies that $\lambda^p=\lambda$ for all $\lambda\in \mathbb F_p$, and it follows that there is a canonical isomorphism $V=\mathbb F\otimes_{\mathbb F_p} V\cong \overline V$ given by 
\[\lambda\otimes v\mapsto \alpha^{-1}(\lambda)\otimes v.\]  
Therefore if $M$ is an $\mathbb F$ vector space, a $p$ semi-linear map $V\to M$ is a linear map $\overline V\to M$ and vice versa.  In this spirit, we use the symbol $\Hom_{\mathbb F}(\overline {\frak g};M)$ to denote the space of all $p$ semi-linear maps of $\frak g$ into $M$.    

The main facts about abelian restricted Lie algebras concern choosing special basis under some additional assumptions.  For us, we will begin our investigations of restricted Lie algebra cohomology in the abelian case since the absence of the Lie bracket in the $p$-operators linearly simplifies the situation considerably.   


\subsection {Restricted Modules}
To study the cohomology of a restricted Lie algebra $\frak g$, it is natural to confine our attention to representations of $\frak g$ that account for the restricted structure.  We continue to fix an arbitrary restricted Lie algebra of characteristic $p>0$.  If $M$ is a vector space over $\mathbb F$, then we will say $M$ is a (restricted) ${\frak g}$-module if there is a restricted Lie algebra homomorphism $\rho:{\frak g}\to \frak {gl}(M)$. Following usual notational conventions, if $g\in \frak g$ and $x\in M$, we will denote the element $\rho(g)(x)\in M$ by $gx$.  We note that the conditions on $\rho$ imply that the pairing $(g,x)\mapsto gx$ is bilinear and that 
\begin{eqnarray*}
[gh]x&=&ghx-hgx \\
g^{[p]}x&=&g^px
\end{eqnarray*}
for all $g,h\in\frak g$ and all $x\in M$. The notions of invariants, irreducibility and complete reducibility for restricted representations are defined precisely as they are for representations of ordinary Lie algebras.

We recall here that if $\frak g$ is an ordinary Lie algebra and $\rho:\frak g\to \frak {gl}(M)$ is an ordinary Lie algebra representation of $\frak g$, then there is a unique associative algebra homomorphism $\tilde\rho:U(\frak g)\to \End M$ satisfying $\tilde\rho\circ \pi=\rho$ where $\pi:\frak g\to U(\frak g)$ is the Poincar\'e-Birkhoff-Witt map.  Moreover, every ordinary (unitary) representation $U(\frak g)\to \End M$ restricts to a Lie algebra representation $\frak g\to \frak {gl}(M)$ so that there is a one-to-one correspondence between representations of the Lie algebra $\frak g$ and (unitary) representations of the universal enveloping algebra $U(\frak g)$.  Finally we recall that this correspondence preserves irreducibility so that in fact, the category of Lie algebra representations of $\frak g$ is naturally isomorphic to the category of unital representations of $U(\frak g)$.  This fact is precisely what allows that ordinary cohomology theory of a Lie algebra $\frak g$ to be defined with a free resolution of the ground field by  $U(\frak g)$-modules.  We now proceed to define the analog of the enveloping algebra for restricted Lie algebras.  That is, we wish to construct an unital associative algebra with an analogous universal mapping property with respect to restricted representations of $\frak g$.  

Let $J$ denote the two-sided ideal in $U(\frak g)$ generated by all elements of the form 
\[g^{[p]}-g^p\]
where $g$ ranges over $\frak g$.  We note that property ({\it R3}) in Definition (\ref{def:maindef}) implies that the generators $g^{[p]}-g^p$ are central in $U(\frak g)$. We denote the quotient $U(\frak g)/J$ by $U_{\rm res.}(\frak g)$ and refer to it as the restricted universal enveloping algebra for the restricted Lie algebra $\frak g$.  We note here that the augmentation $\epsilon:U(\frak g)\to \mathbb F$ vanishes on $J$ and hence induces an augmentation $U_{\rm res.}(\frak g)\to \mathbb F$ which we continue to denote by $\epsilon$.  As for ordinary universal enveloping algebras, we denote the augmentation ideal $\Ker\epsilon$ by $U_{\rm res.}(\frak g)^+$.   We remark that $\mathbb F$ is a trivial $U_{\rm res.}(\frak g)$-module via the action $g\lambda=\epsilon(g)\cdot \lambda$.  We summarize the main properties of the algebra $U_{\rm res.}(\frak g)$ in the following theorem.  The proof(s) can be found in \cite {Jb}, pp. 185 - 192.
\begin{vanish}
\label{thm:ualg}
If $\frak g$ is a restricted Lie algebra of characteristic $p>0$ and $\{e_i\}_{i\in\Lambda}$ is a possibly infinite ordered basis for $\frak g$, then
\begin{itemize}
\item[(1)] The monomials $e_1^{k_1}e_2^{k_2}\cdots e_l^{k_l}$ with $l\ge 0$ and $0\le k_j <p$ form a basis for $U_{\rm res.}(\frak g)$.  In particular, if $\dim_{\mathbb F}\frak g=n$, then $\dim_{\mathbb F}U_{\rm res.}(\frak g)=p^n$.
\item[(2)] The composition $\frak g\stackrel{\scriptscriptstyle\mathrm{PBW}}{\longrightarrow}U(\frak g)\longrightarrow U_{\rm res.}(\frak g)$ is injective.
\item[(3)] There is a one-to-one correspondence between restricted Lie algebra representations of $\frak g$ and unitary representations of $U_{\rm res.}(\frak g)$.  Moreover, this correspondence preserves irreducibility.
\end{itemize}
\end{vanish}
Statement (3) of Theorem (\ref {thm:ualg}) is what we need to parallel the Cartan-Eilenberg construction of ordinary Lie algebra cohomology in the restricted case.  That is, it implies that the cohomology theory of a restricted Lie algebra can be entirely constructed using the associative algebra $U_{\rm res.}(\frak g)$.  We also remark here that together, conditions (1) and (3) of Theorem (\ref {thm:ualg}) imply that every finite dimensional restricted Lie algebra has a finite dimensional faithful representation.  We close this subsection by remarking that the quotient map $U(\frak g)\to U_{\rm res.}(\frak g)$ makes every $U_{\rm res.}(\frak g)$ module a $U(\frak g)$-module so that any resolution of $\mathbb F$ by $U_{\rm res.}(\frak g)$-modules is also a resolution of $U(\frak g)$-modules.


\newpage
\pagestyle{myheadings} 
\markright{  \rm \normalsize CHAPTER 3. \hspace{0.5cm} 
  COHOMOLOGY OF RESTRICTED LIE ALGEBRAS}
\chapter{Restricted Lie Algebra Cohomology}
\thispagestyle{myheadings}
\section {The Abelian Case}
In this chapter, we develop the cohomology theory of restricted Lie algebras, and therefore we consider only restricted Lie algebra modules for coefficients.  As we will see, the presence of the Lie bracket terms in $(g+h)^{[p]}$ complicates the situation considerably so that we begin by considering the abelian case.  Suppose then that we are given an abelian restricted Lie Algebra ${\frak g}$ and a restricted ${\frak g}$-module $M$.  We  will begin with the Cartan-Eilenberg definition of the cohomology groups of ${\frak g}$.  That is, we will attempt to construct a free resolution of $\mathbb F$ by $U_{\rm res.}({\frak g})$-modules where, as before,  $U_{\rm res.}({\frak g})$ denotes the restricted universal enveloping algebra and $\mathbb F$ is regarded as a trivial $U_{\rm res.}(\frak g)$-module.  Of course the Cartan-Eilenberg theory shows that the resulting cohomology groups are independent of the particular resolution that we construct.  In an effort to keep our notations as simple as possible, we will eventually use the notation $C^q({\frak g};M)$ to denote the space of $q$-dimensional cochains of a restricted Lie algebra $\frak g$ with coefficients in a restricted module $M$.  To avoid ambiguity, we will use the notation $-_{\rm cl.}$ when dealing with the cohomology of ${\frak g}$ considered as an ordinary Lie algebra.  For example, the space of  ordinary $q$-dimensional cochains from section (2.1.1) is  $C_{\rm cl.}^q({\frak g};M)$, the ordinary coboundary operator is $\delta_{\rm cl.}$, and so on.  

As we noted in Section (2.2.2), the correspondence between restricted $\frak g$-modules and $U_{\rm res.}(\frak g)$-modules allows us to construct the cohomology theory of a restricted Lie algebra $\frak g$ using the associative algebra $U_{\rm res.}(\frak g)$.  Therefore we proceed to try to construct a free resolution of $\mathbb F$ in the category of $U_{\rm res.}(\frak g)$-modules.  To begin, we set $C_0=U_{\rm res.}(\frak g)$ regarded as a regular $U_{\rm res.}(\frak g)$-module and we have the augmentation $\epsilon:U_{\rm res.}(\frak g)\to \mathbb F$.  If $t,s\ge0$, but not both zero,  we define 
\[ C_{t,s}=S^t\overline{\frak g}\otimes \Lambda^s\frak g\otimes U_{\rm res.}(\frak g)
\]
with its natural structure of a $U_{\rm res.}(\frak g)$-module, and we set 
\[C_k=\bigoplus_{2t+s=k} C_{t,s}.\]
Evidently each $C_k$ is a free $U_{\rm res.}(\frak g)$-module.  We then define a map 
\[d_{t,s}:C_{t,s}\to C_{t,s-1}\oplus C_{t-1,s+1}\]
by the formula
\begin{eqnarray}
\lefteqn{d_{t,s}(h_1\cdots h_t\otimes g_1\wedge\cdots\wedge g_s\otimes x)=} \nonumber \\
& &\sum_{i=1}^s(-1)^{i-1}h_1\cdots h_t\otimes g_1\wedge\cdots\widehat {g_i}\cdots\wedge g_s\otimes g_ix \label {eq:typeone}\\
&+ &  \sum_{j=1}^t h_1\cdots\widehat {h_j}\cdots h_t\otimes h_j^{[p]}\wedge g_1\wedge\cdots\wedge g_s\otimes x \label{eq:typetwo}\\
&- & \sum_{j=1}^t h_1\cdots\widehat {h_j}\cdots h_t\otimes h_j\wedge g_1\wedge\cdots\wedge g_s\otimes h_j^{p-1}x.\label{eq:typethree}
\end{eqnarray}
Clearly each $d_{t,s}$ is a $U_{\rm res.}(\frak g)$-module homomorphism so that the differential 
\[d_k:C_k\to C_{k-1}\]
defined by $d_k=\oplus d_{t,s}$ is also a $U_{\rm res.}(\frak g)$-module map.  
The justification of the term ``differential'' is established by the following theorem.
\begin{theorem}
\label{th:abeliancomplex}
The maps $d_k$ defined above satisfy $d_{k-1}d_k=0$ for $k\ge 1$ and $\epsilon d_1=0$ so that $\{C_k,d_k\}$ is an augmented complex of free $U_{\rm res.}(\frak g)$-modules.
\end{theorem}

\nid {\bf Proof.} The terms in the sum (\ref{eq:typeone}) are elements of $C_{t,s-1}$ whereas the terms in the sums (\ref{eq:typetwo}) and (\ref{eq:typethree}) lie in $C_{t-1,s+1}$.  Therefore, in order to compute $d_{k-1}d_k$, we must apply $d_{t,s-1}$ to (\ref{eq:typeone}) and $d_{t-1,s+1}$ to (\ref{eq:typetwo}) and (\ref {eq:typethree}). Applying $d_{t,s-1}$ to (\ref {eq:typeone}), we have
\begin{eqnarray}
\lefteqn{d_{t,s-1}\left (\sum_{i=1}^s(-1)^{i-1}h_1\cdots h_t\otimes g_1\wedge\cdots\widehat {g_i}\cdots\wedge g_s\otimes g_ix\right ) = }\nonumber \\
& & \sum_{i=1}^{s}(-1)^{i-1}\left (\sum_{\sigma <i}(-1)^{\sigma-1}h_1\cdots h_t\otimes g_1\wedge\cdots\widehat {g_\sigma}\cdots\widehat {g_i}\cdots\wedge g_{s}\otimes g_\sigma g_ix\right. \nonumber \\
&+ & \sum_{\sigma >i}(-1)^{\sigma}h_1\cdots h_t\otimes g_1\wedge\cdots\widehat {g_i}\cdots\widehat {g_{\sigma}}\cdots\wedge g_{s}\otimes g_{\sigma}g_i x \nonumber \\
&+ & \sum_{j=1}^t h_1\cdots\widehat {h_j}\cdots h_t\otimes h_j^{[p]}\wedge g_1\wedge\cdots\widehat {g_i}\cdots\wedge g_{s}\otimes g_i x \label{eq:a} \\
&- & \left.\sum_{j=1}^t h_1\cdots\widehat {h_j}\cdots h_t\otimes h_{j}\wedge g_1\wedge\cdots\widehat {g_i}\cdots\wedge g_{s}\otimes h_j^{p-1}g_ix\right ).\label{eq:b}
\end{eqnarray}
Since $\frak g$ is abelian, $g_\sigma g_i-g_i g_\sigma=0$ in $U_{\rm res.}(\frak g)$ so that the terms in the first two sums in the parentheses cancel in pairs when summed over all $i$.  This leaves the sum over $i$ of (\ref {eq:a}) and (\ref {eq:b}).  Now we apply $d_{t-1,s+1}$ to (\ref{eq:typetwo}).
\begin{eqnarray}
\lefteqn {d_{t-1,s+1}\left (\sum_{j=1}^t h_1\cdots\widehat {h_j}\cdots h_t\otimes h_j^{[p]}\wedge g_1\wedge\cdots\wedge g_s\otimes x    \right ) =}\nonumber \\
\lefteqn {\sum_{j=1}^t \left ( \sum_{\sigma=1}^{s}(-1)^{\sigma} h_1\cdots\widehat {h_j}\cdots h_t\otimes h_j^{[p]}\wedge g_1\wedge\cdots \widehat {g_\sigma}\cdots\wedge g_{s}\otimes g_\sigma x\right. }\label{eq:a2} \\
&+& h_1\cdots\widehat {h_j}\cdots h_t\otimes g_1\wedge\cdots \wedge g_{s}\otimes h_j^{[p]} x \label{eq:c} \\
&+& \sum_{\tau\ne j} h_1\cdots\widehat {h_\tau}\cdots\widehat {h_j}\cdots h_t\otimes h_\tau^{[p]}\wedge h_j^{[p]}\wedge g_1\wedge\cdots\wedge g_{s}\otimes x \label{eq:self1} \\
&-&\left. \sum_{\tau\ne j} h_1\cdots\widehat {h_\tau}\cdots\widehat {h_j}\cdots h_t\otimes h_\tau\wedge h_j^{[p]}\wedge g_1\wedge\cdots\wedge g_{s}\otimes h_\tau^{p-1}x \right ).\label {eq:d} 
\end{eqnarray}
We note that the terms in (\ref {eq:self1}) cancel in pairs since interchanging the first two terms in the alternating product multiplies the term by $-1$.  Finally, we apply $d_{t-1,s+1}$ to (\ref {eq:typethree}) to get
\begin{eqnarray}
\lefteqn {d_{t-1,s+1}\left (- \sum_{j=1}^t h_1\cdots\widehat {h_j}\cdots h_t\otimes h_j\wedge g_1\wedge\cdots\wedge g_s\otimes h_j^{p-1}x \right ) =}\nonumber \\
\lefteqn {- \sum_{j=1}^t\left ( \sum_{\sigma=1}^{s}(-1)^{\sigma} h_1\cdots\widehat {h_j}\cdots h_t\otimes h_j\wedge g_1\wedge\cdots \widehat {g_\sigma}\cdots\wedge g_{s}\otimes g_\sigma h_j^{p-1}x\right. }\label{eq:b2} \\
&+&h_1\cdots\widehat {h_j}\cdots h_t\otimes g_1\wedge\cdots \wedge g_{s}\otimes h_j^{p} x \label{eq:c2} \\
&+& \sum_{\tau\ne j} h_1\cdots\widehat {h_\tau}\cdots\widehat {h_j}\cdots h_t\otimes h_\tau^{[p]}\wedge h_j\wedge g_1\wedge\cdots\wedge g_{s}\otimes h_j^{p-1}x \label{eq:d2} \\ 
&-&\left. \sum_{\tau\ne j} h_1\cdots\widehat {h_\tau}\cdots\widehat {h_j}\cdots h_t\otimes h_\tau\wedge h_j\wedge g_1\wedge\cdots\wedge g_{s}\otimes h_\tau^{p-1}h_j^{p-1}x  \right ).  \label{eq:self2} 
\end{eqnarray}
This time the terms in (\ref {eq:self2}) cancel in pairs.  Moreover, the terms in (\ref {eq:a}) and (\ref {eq:a2}) are identical (with $\sigma=i$) except for sign and hence they cancel.  The terms in (\ref {eq:b}) and (\ref {eq:b2}) cancel in pairs since $h_j^{p-1}g_i-g_i h_j^{p-1}=0$.  The terms in (\ref{eq:d}) and (\ref{eq:d2}) have the same sign but are equal apart from interchanging the first two terms in the alternating part.  Finally the terms in (\ref {eq:c}) and (\ref {eq:c2}) match except for sign since $h_j^{[p]}=h^p_j$ in $U_{\rm res.} (\frak g)$ and hence the entire sum is zero as claimed.  It remains to show that $\epsilon d_1=0$.  But by definition, $C_1=\mathbb F\otimes \frak g\otimes U_{\rm res.}(\frak g)=\frak g\otimes U_{\rm res.}(\frak g)$ and 
\[\epsilon (d_1 (g\otimes x))=\epsilon (gx)=0.\]
This completes the proof.
\qed

We remark here that we may replace $U_{\rm res.}(\frak g)$ in the definition of $C_k$ with any restricted $\frak g$-module $M$ and the above proof goes through unchanged to show that the resulting sequence is a complex.  However, in the case where $M=U_{\rm res.}(\frak g)$, we have the following stronger result.
\begin{theorem}
\label{th:abelianresolution}
If $C=\{C_k,d_k\}$ is the complex defined above, then $H_k(C)=0$ for $0\le k <p$ so that if we define
\[H^k(\frak g;M)=\frac{\Ker(\Hom_{U_{\rm res.}(\frak g)}(C_k,M)\to\Hom_{U_{\rm res.}(\frak g)}(C_{k+1},M))}{\Im(\Hom_{U_{\rm res.}(\frak g)}(C_{k-1},M)\to\Hom_{U_{\rm res.}(\frak g)}(C_{k},M))},\]
then 
\[H^k({\frak g};M)=\Ext^k_{U_{\rm res.}({\scriptstyle \frak g})}(\mathbb F,M)\]
for $0\le k< p$.
\end{theorem}

Our proof of Theorem (\ref {th:abelianresolution}) will require the computation of the homology of two auxiliary complexes.  The first complex is defined as follows.  For each $k\ge 0$, we let $\C_k=C_{0,k}=\Lambda^k\frak g\otimes U_{\rm res.}(\frak g)$ and we define $\partial_{\C}=\partial:\C_k\to \C_{k-1}$ by the formula
\[\partial (g_1\wedge\cdots\wedge g_k\otimes x)=\sum_{i=1}^k (-1)^{i-1}g_1\wedge\cdots\widehat {g_i}\cdots\wedge g_{k}\otimes g_ix.\]
The computation in the proof of Theorem (\ref {th:abeliancomplex}) shows $\partial^2=0$ so that $\{\C_*,\partial\}$ is a complex.  We fix a basis $\{e_1,\dots,e_n\}$ in $\frak g$ once and for all and assume temporarily that $\frak g^{[p]}=0$ so that in particular, $g^p=0$ in $U=U_{\rm res.}(\frak g)$ for all $g\in\frak g$.  If we choose a sequence $1\le i_1<\cdots <i_k\le n$, then it is easy to see that the element 
\begin{equation}
\label{special}
e_{i_1}\wedge\cdots\wedge e_{i_k}\otimes e_{i_1}^{p-1}\cdots e_{i_k}^{p-1}
\end{equation}
is a cycle in $\C_k$ since $e_{i_j}^p=0$ in $U$.  We denote the homology class of this element by $h_{i_1\cdots i_k}$.  In these notations, we have the following lemma.

\begin{lemma}
\label{auxiliarycomplex}
$H_k(\C)$ is spanned by the homology classes
\[\{h_{i_1\cdots i_k}: 1\le i_1<\cdots <i_k\le n\}.\]
\end{lemma}

{\bf Proof.}  For each $k\ge 0$, and each $1\le s\le n$, we let $\C_{k,s}$ be the $\mathbb F$-subspace of $\C_k$ spanned by all elements of the form
\begin{equation}
\label{subcomplexes}
e_{i_1}\wedge\cdots\wedge e_{i_l}\wedge e_{j_1}\wedge\cdots\wedge e_{j_m}\otimes e_{i_1}^{p-1}\cdots e_{i_l}^{p-1}e_{s+1}^{r_{s+1}}\cdots e_n^{r_n}
\end{equation}
where 
\[\begin{array}{c}
1\le i_1<\cdots < i_l\le s <j_1<\cdots <j_m\le n\\
l+m=k\\
0\le r_a\le p-1\ \hbox{\rm for}\ s+1\le a \le n.
\end{array}
\]
Clearly the boundary of an element of the form (\ref {subcomplexes}) is either 0 or has the same form so that we have a nested sequence of subcomplexes
\[\C=\C_{\cdot 0}\supset\C_{\cdot 1}\supset\cdots\supset\C_{\cdot n}\]
where $\C_{\cdot s}=\bigoplus_k \C_{k,s}$.
We claim that the inclusion map $\C_{\cdot s}\to \C_{\cdot s-1}$ induces an isomorphism in homology for all $s$, and hence $H(\C)=H(\C_{\cdot n})$. To see this, we introduce a filtration in the quotient complex $\C_{\cdot s-1}/\C_{\cdot s}$ as follows.  We define $F_t=F_t(\C_{\cdot s-1}/\C_{\cdot s})$ to be the subspace of $\C_{\cdot s-1}/\C_{\cdot s}$ spanned by monomials in which the total degree of the exterior part in $e_1,\dots,\widehat {e_s},\dots, e_n$ is less than or equal to $t$.  Since the boundary operator $\partial$ either preserves this degree or decreases it by one, we see that this filtration is compatible with the differential $\partial$ and hence we have a corresponding (homology) spectral sequence with $E^0_t=F_t/F_{t-1}$.  If we denote the induced differential $E^0_t\to E^0_t$ by $\partial_s$, then it is easy to see that $\partial_s$ acts on monomials $a\otimes u$ by the formula
\begin{equation*}
\partial_s(a\otimes u)=\left \{
\begin{array}{cl}
b\otimes e_su& \hbox {\rm if $a$ contains $e_s$;}\\
0 & \hbox {\rm if $a$ does not contain $e_s$}.
\end{array}\right.
\end{equation*}
where we write $a=e_s\wedge b$; that is $b\in\Lambda^*\frak g$ is the monomial that results from deleting $e_s$ from $a$.  Let $\deg_{e_s} (a\otimes u)$ denote the total degree of $e_s$ in the monomials $u\in U, a\in \Lambda^*\frak g$ and note that if $a\otimes u$ represents an element of the quotient $\C_{\cdot s-1}/\C_{\cdot s}$, then $\deg_{e_s}(a\otimes u)\ne 0,p$. We then define a map $D_s:E^0_t\to E^0_t$ on monomials $a\in \Lambda^*\frak g$ and $u\in U$ by the formula
\begin{equation*}
D_s(a\otimes u)=a\wedge e_s\otimes \left (\frac{1}{\deg_{e_s}(a\otimes u)}\right )\frac{\partial u}{\partial e_s}
\end{equation*}
where  and $\partial u/\partial e_s$ denotes the derivative of the monomial $u$ with respect to the variable $e_s$. If we set $\deg_{e_s} u=\deg_{e_s} (1\otimes u)$, then we have the equality
\[e_s\frac{\partial u}{\partial e_s}=(\deg_{e_s} u)u.\]
Now we compute for any tensor product of monomials $a\otimes u\in E^0_t$, 
\begin{equation*}
\partial_s D_s(a\otimes u)=\left \{
\begin{array}{cl}
a\otimes \frac{1}{\deg_{e_s}(a\otimes u)}\deg_{e_s}(u)u    & \hbox {\rm if $a$ does not contain $e_s$;}\\
0 & \hbox {\rm if $a$ contains $e_s$}.
\end{array}\right.
\end{equation*}
But if the monomial $a$ does not contain $e_s$, then $\deg_{e_s}(a\otimes u)=\deg_{e_s}u$ so that we have 
\begin{equation}
\label{littlebig}
\partial_s D_s(a\otimes u)=\left \{
\begin{array}{cl}
a\otimes u  & \hbox {\rm if $a$ does not contain $e_s$;}\\
0 & \hbox {\rm if $a$ contains $e_s$}.
\end{array}\right.
\end{equation}
On the other hand, we have
\begin{equation*}
 D_s\partial_s(a\otimes u)=\left \{
\begin{array}{cl}
0& \hbox {\rm if $a$ does not contain $e_s$;}\\
\frac{1}{\deg_{e_s}(b\otimes e_su)}\left (a\otimes u+a\otimes \deg_{e_s}(u)u\right )& \hbox {\rm if $a$ contains $e_s$}
\end{array}\right.
\end{equation*}
where we write $a=e_s\wedge b$ as before.  Now, $\deg_{e_s}(b\otimes e_su)=1+\deg_{e_s}(u)$ so that we have 
\begin{equation}
\label{biglittle}
 D_s\partial_s(a\otimes u)=\left \{
\begin{array}{cl}
0& \hbox {\rm if $a$ does not contain $e_s$;}\\
a\otimes u& \hbox {\rm if $a$ contains $e_s$}.
\end{array}\right.
\end{equation}
Therefore if we add the results from equations (\ref {littlebig}) and (\ref {biglittle}), we have 
\[\partial_sD_s+D_s\partial_s=1_{E^0_t}\]
so that $E^1_t=E^\infty_t =0$ and the spectral sequence of the filtration $F$ collapses completely.  It follows immediately the quotient complex $\C_{\cdot s-1}/\C_{\cdot s}$ is acyclic and hence the homology $H(\C_{\cdot s-1})$ is equal to the homology $H(\C_{\cdot s})$ for all $s$. In particular, we have $H(\C)=H(\C_{\cdot n})$.  It is obvious from the definitions that the boundary operator on $\C_{\cdot n}$ is identically zero so that $H_k(\C_{\cdot n})=\C_{k,n}$ and the latter is spanned by cycles of the form (\ref {special}) by definition.  This completes the proof of the lemma.
\qed

We remark that our proof of Lemma (\ref {auxiliarycomplex}) actually shows that the homology classes $h_{i_1\cdots i_k}$ for a basis form $H(\C)$ so that $\dim H_k(\C)=\bino{n}{k}$.  In the case $\frak g^{[p]}\ne 0$, the elements (\ref {special}) are not necessarily cycles.  However, the elements
\[c_i=e_i^{[p]}\otimes 1-e_i\otimes e_i^{p-1}\in \C_1\] are clearly cycles and if we define the product of such elements with the usual product in the tensor product of $\mathbb F$-algebras, we have 
\begin{equation}
\label{sospecial}
c_{i_1}\cdots c_{i_k}=\sum_{J\subset \{1,\dots ,k\}}(-1)^{|J|} f_{i_1}\wedge\dots \wedge f_{i_k}\otimes e_{i_1}^{q_{i_1}}\cdots e_{i_k}^{q_{i_k}}
\end{equation}
where $1\le i_1<\cdots < i_k\le n$ and 
\begin{eqnarray*}
f_{i_j}=\left \{\begin{array}{cl} e_{i_j}& \hbox {\rm if $j\in J$}\\ e_{i_j}^{[p]} & \hbox {\rm if $j\notin J$}\end{array}\right. & \hbox {\rm and} & q_{i_j}=\left \{\begin{array}{cl} p-1 & \hbox {\rm if $j\in J$}\\ 0&  \hbox {\rm if $j\notin J$.}\end{array}\right.
\end{eqnarray*}
Evidently each $c_{i_1}\cdots c_{i_k}$ is a cycle and we denote its homology class by $\tilde h_{i_1\cdots i_k}$.  We then let $\tilde\C_{k,n}$ be the $\mathbb F$-subspace of $\C_k$ spanned by the elements $c_{i_1}\cdots c_{i_k}$ where $1\le i_1<\cdots <i_k\le n$ and $\tilde\C_{\cdot n}=\bigoplus_k \tilde\C_{k,n}$.  We claim that the inclusion map $\tilde\C_{\cdot n}\to \C$ induces an isomorphism in homology and hence $H(\C)$ is spanned by the homology classes $\tilde h_{i_1\cdots i_k}$.  A cycle $c\in \C$ is a sum of monomials
\begin{equation}
\label{monomial}
e_I\otimes e^r=e_{i_1}\wedge \cdots\wedge e_{i_k}\otimes e_1^{r_1}\dots e_n^{r_n}\end{equation}
where $k\ge 0$, $I=(i_1,\dots,i_k)$ is increasing and $r=(r_1,\dots,r_n)$ satisfies $0\le r_j\le p-1$.  We define the total degree of a monomial $e_I\otimes e^r$ by 
\[\Deg(e_I\otimes e^r)=k+\sum_j r_j.\]
We remark that the boundary operator $\partial:\C\to\C$ either preserves $\Deg(e_I\otimes e^r)$ or lowers it by $p-1$ so that we can write $\partial=\partial^0+\partial^-$ where $\partial^0$ preserves the total degree and $\partial^-$ lowers it by $p-1$.  We note that $\partial^0$ is precisely the boundary operator on the complex $\C_{\cdot 0}$.  To show that $H(\tilde\C_{\cdot n})\to H(\C)$ is an epimorphism, it suffices to show that given a cycle $c\in C$, we can find a cycle $\tilde c\in \tilde\C_{\cdot n}$ such that $c-\tilde c=\partial b$ for some chain $b\in \C$.  Given a chain $c\in\C$, let $c^*$ denote the sum of the monomials in $c$ of maximal total degree.  We note that if $\partial c=0$, then $\partial^0 c^*=0$ since $\partial^0$ preserves total degree.  Consequently, Lemma (\ref {auxiliarycomplex}) implies that $c^*-c'=\partial^0 b$ where $c'$ is a cycle spanned by monomials of the form (\ref {special}) and $\Deg (b)=\Deg (c)$.  It follows that $(c-\partial b)^*$ is a cycle in $\tilde\C_{\cdot n}$, say $(c-\partial b)^*=\tilde c$.  If we let $c''=c-\partial b- \tilde c$, then we have $\Deg (c'')<\Deg (c)$ and $c-\tilde c=\partial b - c''$.  Inducting on $\Deg c$ then shows that $c-\tilde c=\partial b$ for some cycle $\tilde c\in\tilde\C_{\cdot n}$ and hence our map is an epimorphism as claimed.   Similarly we can show that this map is a monomorphism so that we have shown the following lemma. 
\begin{lemma}
\label{notzerocase}
If $\C$ is the complex defined above, then $H(\C)=H(\tilde\C_{\cdot n})$ so that $H_k(\C)$ has a basis consisting of the homology classes 
\[\{\tilde h_{i_1\cdots i_k} : 1\le i_1<\cdots < i_k\le n\}.\]
In particular we have $\dim H_k(\C)=\bino{n}{k}$.
\qed
\end{lemma} 

A basis for the space $C_{t,s}$ consists of the monomials
\begin{equation}
\label{genmonomial}
e^\mu\otimes e_I\otimes e^r= e_1^{\mu_1}\cdots e_n^{\mu_n}\otimes e_{i_1}\wedge \cdots \wedge e_{i_s}\otimes e_1^{r_1}\cdots e_n^{r_n}
\end{equation}
where $\mu=(\mu_1,\dots,\mu_n)$, $I=(i_1,\dots,i_s)$, $r=(r_1,\dots, r_n)$ and
\[\begin{array}{c} 0\le \mu_j\\
|\mu|=\sum_j\mu_j=t\\
1\le i_1<\cdots <i_s\le n\\
0\le r_j\le p-1.
\end{array} 
\]
For each $i=1,\dots, n$, we let
\[c_i=1\otimes e_i^{[p]}\otimes 1-1\otimes e_i\otimes e_i^{p-1}\]
 and we easily note that $c_i\in C_{0,1}$ is a cycle for all $i$.  Now we define
\[(\partial/\partial e_i\otimes c_i):C_{t,s}\to C_{t-1,s+1}\]
 by the formula
\[\left (\frac{\partial}{\partial e_i}\otimes c_i\right )(e^\mu\otimes e_I\otimes e^r)=\frac{\partial e^\mu}{\partial e_i}\otimes e_i^{[p]}\wedge e_I\otimes e^r-\frac{\partial e^\mu}{\partial e_i}\otimes e_i\wedge e_I\otimes e_i^{p-1}e^r.\]
If $\mu=(\mu_1,\dots,\mu_n)$ satisfies $|\mu|=t$ and $I=(i_1,\dots,i_s)$ is increasing, then by definition we write
\[e^\mu\otimes c_I=\sum_{J\subset \{1,\dots ,s\}}(-1)^{|J|} e^\mu\otimes f_{i_1}\wedge\dots \wedge f_{i_s}\otimes e_{i_1}^{q_{i_1}}\cdots e_{i_s}^{q_{i_s}}\]
where 
\begin{eqnarray*}
f_{i_j}=\left \{\begin{array}{cl} e_{i_j}& \hbox {\rm if $j\in J$}\\ e_{i_j}^{[p]} & \hbox {\rm if $j\notin J$}\end{array}\right. & \hbox {\rm and} & q_{i_j}=\left \{\begin{array}{cl} p-1 & \hbox {\rm if $j\in J$}\\ 0&  \hbox {\rm if $j\notin J$.}\end{array}\right.
\end{eqnarray*}
We then define $\frak C_{t,s}$ to be the $\mathbb F$-subspace of $C_{t,s}$ spanned by the elements
\[\{e^\mu\otimes c_I : \hbox {\rm $|\mu|=t$ and $I$ is increasing}\}\]
and 
\[\frak C_k=\bigoplus_{2t+s=k} \frak C_{t,s}.\]
The boundary operator $\partial_{\frak C}=\partial:\frak C_k\to \frak C_{k-1}$ is defined by
\[\partial=\sum_{j=1}^n \frac{\partial}{\partial e_j}\otimes c_j.\]
The proof of Theorem (\ref {th:abeliancomplex}) shows that $\partial^2=0$ so that $\frak C=\{\frak C_*,\partial\}$ is a complex.  

\begin{lemma}
\label{auxiliarytwo}
If $\frak C$ is the complex defined above, then 
\[H_k(\frak C)=\left\{\begin{array}{cl} U_{\rm res.}(\frak g)& \hbox {\rm if $k=0$}; \\
0 & \hbox {\rm if $0<k<p$}
\end{array}\right.\]
\end{lemma}

{\bf Proof.} Define a map $D:\frak C_k\to \frak C_{k+1}$ by the formula
\[D(e^\mu\otimes c_I)=\sum_{a=1}^s (-1)^{a-1} e^\mu e_{i_a}\otimes c_{i_1}\cdots \widehat {c_{i_a}}\cdots c_{i_s}\]
and compute for any monomial $e^\mu\otimes c_I$:
\begin{eqnarray}
\lefteqn{D\partial(e^\mu\otimes c_I) = \sum_{\stackrel{\scriptstyle j=1}{j\ne i_1,\dots,i_s}}^n D(\mu_j e_1^{\mu_1}\cdots e_j^{\mu_j-1}\cdots e_n^{\mu_n}\otimes c_jc_I)}\nonumber \\
&=&\left (\sum_{\stackrel{\scriptstyle j=1}{j\ne i_1,\dots,i_s}}^n\mu_j\right )e^\mu\otimes c_I \label{degrees1}\\
&-& \sum_{\stackrel{\scriptstyle j=1}{j\ne i_1,\dots,i_s}}^n\sum_{a=1}^s(-1)^{a}\mu_j e_1^{\mu_1}\cdots e_j^{\mu_j-1}\cdots e_{i_a}^{\mu_{i_a}+1}\cdots e_n^{\mu_n}\otimes c_jc_{i_1}\cdot \widehat {c_{i_a}}\cdot c_{i_s}\label{alpha}
\end{eqnarray}
and 
\begin{eqnarray}
\lefteqn{\partial D(e^\mu\otimes c_I) = \sum_{a=1}^s (-1)^{a-1}\partial(e_1^{\mu_1}\cdots e_{i_a}^{\mu_{i_a}+1}\cdots e_n^{\mu_n}\otimes c_{i_1}\cdots \widehat {c_{i_a}}\cdots c_{i_s})}\nonumber \\
&=&\left (\sum_{a=1}^s \mu_{i_a}+1\right )e^\mu\otimes c_I \label{degrees2}\\
&+& \sum_{a=1}^s(-1)^{a}\sum_{\stackrel{\scriptstyle j=1}{j\ne i_1,\dots,i_s}}^n\mu_j e_1^{\mu_1}\cdots e_j^{\mu_j-1}\cdots e_{i_a}^{\mu_{i_a}+1}\cdots e_n^{\mu_n}\otimes c_jc_{i_1}\cdot \widehat {c_{i_a}}\cdot c_{i_s}.\label{beta}
\end{eqnarray}
Clearly the terms (\ref {alpha}) and (\ref {beta}) are identical apart from sign so that we have
\begin{equation*}
(D\partial+\partial D)(e^\mu\otimes c_I)=\left (\sum_{\stackrel{\scriptstyle j=1}{j\ne i_1,\dots,i_s}}^n\mu_j+\sum_{a=1}^s \mu_{i_a}+s\right )(e^\mu\otimes c_I)=(t+s)(e^\mu\otimes c_I).
\end{equation*}
Therefore we see that every cycle in $\frak C_k$ ($k=2t+s$) is a boundary provided that $s+t\ne 0\pmod {p}$.  In particular, if $0<k<p$, then $0<t+s<p$ so that $H_k(\frak C)=0$.  Moreover, $\frak C_1=\frak C_{0,1}$ is spanned by the $c_i$ and $\partial c_i=0$ for all $i$.  Therefore $H_0(\frak C)=\frak C_0=U_{\rm res.}(\frak g)$, and the proof of the lemma is complete.
\qed

{\bf Proof of Theorem (\ref {th:abelianresolution}).}  Let $c\in C_k$ be a cycle so that $dc=0$.  We know $c$ is a sum of monomials
\[e^\mu\otimes e_I\otimes e^r=e^\mu\otimes A_\mu\]
where $A_\mu\in\C_s$ with  $s=|I|=k-2|\mu|$.  Let $t^*$ be the maximum value of $|\mu|$ in the monomials that comprise $c$ and write
\[c=\sum_{|\mu|=t^*}e^\mu\otimes A_\mu+\sum_{|\mu|<t^*}e^\mu\otimes A_\mu.\]
We claim that $\partial_{\C}A_\mu=0$ if $|\mu|=t^*$.  To see this, notice that we can write the boundary operator $d$ as a sum $d=\partial_{\C}+\partial_{\frak C}$ and $|\mu|$ is invariant with respect to  $\partial_{\C}$ whereas $\partial_{\frak C}$ decreases $|\mu|$ by 1.  Therefore we have 
\[0=dc=\sum_{|\mu|=t^*}e^\mu\otimes \partial_{\C}A_\mu+(\hbox {\rm terms with $|\mu|<t^*$})\]
and it follows that $\partial_{\C}A_\mu=0$ as claimed.  Now, for each $|\mu|=t^*$, Lemma (\ref {notzerocase}) implies that $A_\mu=B_\mu+\partial_{\C}C_\mu$ where $B_\mu$ is a linear combination of cycles of the form (\ref {sospecial}).  Therefore we have 
\[d(e^\mu\otimes C_\mu)=e^\mu\otimes (A_\mu-B_\mu) +(\hbox {\rm terms with $|\mu|<t^*$})\]
so that 
\[e^\mu\otimes A_\mu-d(e^\mu\otimes C_\mu)=e^\mu\otimes B_\mu+(\hbox {\rm terms with $|\mu|<t^*$}).\]
We remark that the terms $e^\mu\otimes B_\mu$ are elements of the complex $\frak C$ by definition.  All of this implies that we lose no generality in assuming the leading term 
\[\sum_{|\mu|=t^*}e^\mu\otimes A_\mu\]
of $c$ is an element of $\frak C$.  We claim that this term is a cycle in $\frak C$. Indeed, we have 
\begin{eqnarray}
\lefteqn{0=dc=\partial_{\C}c+\partial_{\frak C}c =\sum_{|\mu|=t^*} e^\mu\otimes \partial_{\C}A_\mu\nonumber} \\
& & +\sum_{|\mu|<t^*} e^\mu\otimes \partial_{\C}A_\mu\label {egy} \\
& & +\sum_{|\mu|=t^*} \partial_{\frak C}(e^\mu\otimes A_\mu) \label {ketto}\\
& & +\sum_{|\mu|<t^*} \partial_{\frak C}(e^\mu\otimes A_\mu).\nonumber 
\end{eqnarray}
We want to show that (\ref {ketto}) is zero.  Note that each term in (\ref {ketto}) has degree $(t^*-1)$ in the symmetric part and the only other terms of this degree can come from (\ref {egy}).  The terms of degree $(t^*-1)$ in (\ref {egy}) have the form $e^\mu\otimes \partial_{\C}A_\mu$ whereas the terms in (\ref {ketto}) are all in $\frak C$.  But an element of $\frak C$ never looks like $e^\mu\otimes \partial_{\C}A_\mu$ so that the sum (\ref {ketto}) vanishes as claimed.  Now, if $0< k< p$, it follows from Lemma (\ref {auxiliarytwo}) that
\[c=\partial_{\frak C}b+\sum_{|\mu|<t^*}(e^\mu\otimes A_\mu).\]
Note that since $b\in\frak C$, $\partial_{\C}b=0$ so that we have 
\[c-db=c-\partial_{\C}b-\partial_{\frak C}b=\sum_{|\mu|<t^*}(e^\mu\otimes A_\mu).\]
Therefore if we let $c'=\sum_{|\mu|<t^*}(e^\mu\otimes A_\mu)$, we see that $d(c')=0$, the degree of the symmetric part of $c'$ is less than that of $c$ and 
\[c=db+c'.\]
Induction on the degree of the symmetric part of $c$ then shows that if $0<k<p$, then $H_k(C)=0$. If $c\in C_0$ is in the kernel of $\epsilon$, then $c$ is a sum of monomials $1\otimes 1\otimes e^r$ with $r=(r_1,\dots, r_n)$ and not all $r_j=0$. But if $r_j\ne 0$, then $1\otimes e_j\otimes e_1^{r_1}\cdots e_j^{r_j-1}\cdots e_n^{r_n}\in C_1$ and clearly 
\[d(1\otimes e_j\otimes e_1^{r_1}\cdots e_j^{r_j-1}\cdots e_n^{r_n})=1\otimes 1\otimes e^r.\]
Therefore $H_0(C)=0$ and the proof is complete.
\qed

Now if $M$ is a restricted $\frak g$ module, and $0\le k <p$, then Theorem (\ref {th:abelianresolution}) states that 
\[H^k({\frak g};M)=\Ext^k_{U_{\rm res.}({\scriptstyle \frak g})}(\mathbb F,M)\]
so that, using the complex $C$ of Theorem (\ref{th:abelianresolution}) we have
\begin{eqnarray*}
\lefteqn{C^k(\frak g;M)=\Hom_{U_{\rm res.}({\scriptstyle \frak g})}(C_k,M)} \\
& & =\bigoplus_{2t+s=k}\Hom_{U_{\rm res.}({\scriptstyle \frak g})}(S^t\overline {\frak g}\otimes \Lambda^{s}\frak g\otimes U_{\rm res.}(\frak g),M) \\
& & =\bigoplus_{2t+s=k} \Hom_{\mathbb F} (S^t\overline {\frak g}\otimes \Lambda^{s}\frak g,M).
\end{eqnarray*} 
Therefore if $\frak g$ and $M$ are finite dimensional,  we have the following important corollary to Theorem (\ref{th:abelianresolution}).
\begin{corollary}
\label{cor:dimensions}
If $\dim_{\mathbb F}(\frak g)=n$ and $\dim_{\mathbb F}(M)=m$, then the dimension of the space of $k$-dimensional cochains of $\frak g$ with coefficients in $M$ is 
\[\bino{n+k-1}{k}\cdot m\]
In particular the dimension of $C^k(\frak g;\mathbb F)$ is the same as that of the symmetric algebra $S^k\frak g$.
\end{corollary}

{\bf Proof.} Recalling that $\dim_{\mathbb F}\Lambda^s\frak g=\bino{n}{s}$ and $\dim_{\mathbb F}S^t\frak g=\bino{n+t-1}{t}$, the above remarks show that 
\[\dim_{\mathbb F}C^k(\frak g;M)=\left ( \sum_{s+2t=k}\bino{n}{s}\bino{n+t-1}{t}\right )\cdot m.\]
The result then follows from the identity
\begin{equation}
\label{eq:homework1}
\sum_{s+2t=k}\bino{n}{s}\bino{n+t-1}{t}=\bino{n+k-1}{k}.
\end{equation}
To see this, we recall the well know identities
\[(1+t)^n=\sum_{k\ge 0}\bino{n}{k}t^k\]
and 
\[(1-t)^{-n}=\sum_{k\ge 0}\bino{n+k-1}{k}t^k.\]
Now, equating coefficients of $t^k$ in the equality
\[(1+t)^n(1-t^2)^{-n}=(1-t)^{-n}\]
gives the identity (\ref {eq:homework1}) and the proof of the corollary is complete.
\qed

In the subsequent section, we will give explicit constructions of the cochain spaces $C^k(\frak g;M)$ in the nonabelian case for $k\le 3$.  We remark on this here to point out that the dimensions of these spaces will agree with those computed in Corollary (\ref{cor:dimensions}).  This is remarkable because, in the characteristic zero case, an arbitrary Lie algebra is a deformation of an abelian Lie algebra.  Upon defining a suitable notion of deformations of restricted Lie algebras, we might expect that the entire cohomology theory in the non-abelian case may be a deformation of the theory constructed above.  That is, we may be able to think of the cohomology of non-abelian restricted Lie algebra as a sort of quantization of the cohomology of abelian Lie algebras.  


\section {The Complex}
We have seen that in the absence of the Lie bracket in the additive property for the $p$-operator,  we can construct a free resolution of the ground field $\mathbb F$ by $U_{\rm res.}(\frak g)$-modules and therefore construct a complex $\{C^*({\frak g};M),\delta\}$ with cohomology 
\[H^k({\frak g};M)=\Ext^k_{U_{\rm res.}({\scriptstyle \frak g})}(\mathbb F,M).\]   The situation for nonabelian restricted Lie algebras is considerably more complicated.   We will therefore content ourselves with explicit constructions of the $k$-dimensional cochain spaces $C^k({\frak g};M)$ and coboundary operators \[ \delta^k:C^k({\frak g};M)\to C^{k+1}({\frak g};M) \] 
for $k\le 2$.   Continuing to use the notations above, we fix a field $\mathbb F$ with $\mathop {\rm char} \mathbb F=p>0$, a restricted Lie algebra ${\frak g}$ and a restricted ${\frak g}$-module $M$. For $k\le 1$, define $C^k({\frak g};M)=C_{\rm cl.}^k({\frak g};M)$ and $\delta^0:C^0({\frak g};M)\to C^1({\frak g};M)$
by $\delta^0=\delta^0_{\rm cl.}$ It follows immediately that $H^0({\frak g};M)=H_{\rm cl.}^0({\frak g};M)$.  If $\phi:\Lambda^2{\frak g}\to M$ is a skew-symmetric bilinear form on ${\frak g}$ with values in $M$ and $\omega:{\frak g}\to M$, we say that $\omega$ has the {\it $*$-property with respect to $\phi$} if for all $\lambda\in \mathbb F$ and all $g,h\in {\frak g}$, 
\begin{itemize}
\item[(\romannumeral 1)] $\omega(\lambda g)=\lambda^p\omega(g)$.
\item[(\romannumeral 2)] \begin{eqnarray*}
\omega(g+h)&=&\omega(g)+\omega(h) \\
&+&\sum_{\stackrel{g_j=g\ {\rm or}\  h}{\scriptscriptstyle g_1=g,g_2=h}} \frac{1}{\#(g)}\sum_{k=0}^{p-2}(-1)^k g_p\cdots g_{p-k+1}\phi([g_1,\dots,g_{p-k-1}],g_{p-k}).
\end{eqnarray*} 
\end{itemize}
Our space of 2-dimensional cochains is then defined as 
\[ C^2({\frak g};M)= \{(\phi,\omega)\ |\ \phi:\Lambda^2{\frak g}\to M, \omega:{\frak g}\to M\ \hbox{\rm has the}\ *\hbox{\rm -property w.r.t}\  \phi\}.\]
Evidently if $\omega$ and $\omega^\prime$ have the $*$-property with respect to  $\phi$ and $\phi^\prime$ respectively, then $\omega+\omega^\prime$ has the $*$-property with respect to  $\phi+\phi^\prime$, and hence $C^2({\frak g},M)$ is a vector space over $\mathbb F$ by point wise addition in each coordinate.  We remark here that given $\phi\in C^2_{\rm cl.}({\frak g};M)$, there are numerous maps $\omega:{\frak g}\to M$ with the $*$-property with respect to $\phi$.  In fact, if we choose a basis in ${\frak g}$, then the values of $\omega$ on the basis can be assigned arbitrarily and conditions (\romannumeral 1) and (\romannumeral 2) above completely determine $\omega$ for a given $\phi$.  Moreover, the mapping $(\phi,\omega)\mapsto \phi$ of $C^2({\frak g};M)$ into $C^2_{\rm cl.}({\frak g};M)$ is clearly linear, and if we temporarily denote its kernel by $K$, we have an exact sequence of $\mathbb F$ vector spaces
\[0\to K\to C^2({\frak g};M)\to C^2_{\rm cl.}({\frak g};M)\to 0.\]
The kernel $K$ consists of pairs $(0,\omega)$, and $\omega:{\frak g}\to M$ has the $*$-property with respect to the zero map if and only if $\omega:\frak g\to M$ is  $p$ semi-linear as defined in (2.2.1).  Using the notation there, we write $K=\Hom_{\mathbb F}(\overline{\frak g},M)$.  It follows that if $\dim_{\mathbb F}{\frak g}=n$ and $\dim_{\mathbb F} M=m$, then 
\[\dim_{\mathbb F}C^2({\frak g};M)=nm+\bino{n}{2}m=\frac{n(n+1)}{2}m\]
which agrees with the result in Corollary (\ref{cor:dimensions}).
To define the coboundary operator $\delta^1:C^1({\frak g};M)\to C^2({\frak g};M)$, we will need the following lemma.

\begin{lemma}
\label{starprop}
Given a linear map $\psi:{\frak g}\to M$, if we define $\tilde\psi:{\frak g}\to M$ by the formula 
\[ \tilde\psi(g)=\psi(g^{[p]})-g^{p-1}\psi(g), \] 
then $\tilde\psi$ has the $*$-property with respect to $\delta_{\rm cl.}\psi$. 
\end {lemma}

The proof of Lemma (\ref{starprop}) will require the following technical lemma whose proof is an easy induction and is therefore omitted.

\begin{lemma}
\label{bracket}
Let $I= \{i_1,\dots, i_t\}$, $J=\{j_1,\dots, j_s\}$, $r\ge 2$ and write $I+J=\{2,\dots,r\}$  if $I\cup J=\{2,\dots, r\}$, $I\cap J=\emptyset$ and 
\[ 2\le i_1<i_2<\cdots<i_t\le r\hskip.5in {\rm  and}\hskip.5in 2\le j_1<i_2<\cdots<j_s\le r.\] 
It's okay for one of $I$ or $J$ to be empty.  In this notation, we have
$$[g_1,\dots,g_r]=\sum_{I+J=\{2,\dots,r\}} (-1)^{s}g_{j_s}\dots g_{j_1}g_1g_{i_1}\dots g_{i_t}.$$
\qed
\end{lemma}

\nid {\bf Proof of Lemma (\ref{starprop}).} By definition and repeated application of the identity
$$\psi[g_1,\dots,g_p]=-\delta_{\rm cl.}\psi([g_1,\dots,g_{p-1}],g_p)-[g_1,\dots,g_{p-1}]\psi(g_p)+g_p\psi[g_1,\dots,g_{p-1}]$$
we have
\begin{eqnarray}
\tilde\psi(g+h)&=&\psi(g^{[p]})+\psi(h^{[p]}) \nonumber \\
 & + &\sum_{\stackrel{g_j=g\ {\rm or}\ h}{\scriptscriptstyle g_1=g,g_2=h}} \frac{1}{\#(g)}\left (\sum_{k=0}^{p-2}(-1)^k g_p\cdots g_{p-k+1}\delta_{\rm cl.}\psi([g_1,\dots,g_{p-k-1}],g_{p-k})\right ) \nonumber \\
 & + &\sum_{\stackrel{g_j=g\ {\rm or}\ h}{\scriptscriptstyle g_1=g,g_2=h}} \frac{1}{\#(g)}\left (\sum_{k=0}^{p-1}(-1)^k g_p\cdots g_{p-k+1}[g_1,\dots,g_{p-k-1}]\psi(g_{p-k})\right ) \nonumber \\
& - & (g+h)^{p-1}(\psi(g)+\psi(h)). \nonumber
\end{eqnarray}
Therefore it suffices to show that 
\begin{eqnarray}
\label{eq:0}
 & &\sum_{\stackrel{g_j=g\ {\rm or}\ h}{\scriptscriptstyle g_1=g,g_2=h}} \frac{1}{\#(g)}\left (\sum_{k=0}^{p-1}(-1)^k g_p\cdots g_{p-k+1}[g_1,\dots,g_{p-k-1}]\psi(g_{p-k})\right ) \\
&-&(g+h)^{p-1}(\psi(g)+\psi(h))\nonumber  \\
 & &= -g^{p-1}\psi(g)-h^{p-1}\psi(h).\nonumber
\end{eqnarray}
Now by Lemma \ref{bracket}, we have
\begin{eqnarray}
\label{eq:1}
 &\displaystyle\sum_{k=0}^{p-1}(-1)^k g_p\cdots g_{p-k+1}[g_1,\dots,g_{p-k-1}]\psi(g_{p-k})\nonumber \\
=& \displaystyle\sum_{\stackrel{A+B}{\scriptscriptstyle =\{2,\dots,p\}}} (-1)^{s}g_{a_s}\dots g_{a_1}g_1g_{b_1}\dots g_{b_{t-1}}\psi(g_{b_t}).
\end{eqnarray}
Moreover, it is evident that 
\begin{eqnarray}
\label{eq:2}
(g+h)^{p-1}(\psi(g)+\psi(h))=\sum_{h_j=g\ {\rm or}\ h} h_1\cdots h_{p-1}\psi(h_p).
\end{eqnarray}
Now we count how many times $h_1\cdots h_{p-1}\psi(h_p)$ occurs in (\ref {eq:1}).  Suppose that $\# (g)\ne 0$ or $p$ so that there is at least one $g$ and $h$ among the $h_i$, say $g=h_{s+1}$ with $0\le s\le p-1$.  For each such $g$, we have
\[ \bino {p-1}{s}\equiv (-1)^s \pmod {p} \]
choices for $a_1, \dots, a_s\in \{2,\dots, p\}$.  Therefore there are $(-1)^s\#(g)\pmod {p}$ occurrences of $h_1\cdots h_{p-1}\psi(h_p)$ in (\ref{eq:1}).  Since each such term has a coefficient of $(-1)^s$, we see that there are $\#(g)\pmod {p}$ occurrences total.  Therefore if $\#(g)\ne 0$ or $p$, we see that the terms in the first two lines of equation (\ref {eq:0}) cancel in pairs. The remaining terms in (\ref{eq:2}) are $-g^{p-1}\psi(g)$ ($\#(g)=p$) and $-h^{p-1}\psi(h)$ ($\#(g)=0$).  This proves the lemma.
\qed

We can now define the coboundary operator
\[ \delta^1:C^1({\frak g};M)\to C^2({\frak g};M) \]
by $\delta^1:\psi\mapsto (\delta_{\rm cl.}\psi,\tilde\psi)$, and we have our first important result.
\begin{theorem} 
\label{thm:1}
In the above notations, $\delta^1\delta^0=0$ and $H^1({\frak g};M)=\Ker \delta^1/\Im\delta^0$ injects as a subspace of $H_{\rm cl.}^1({\frak g};M)$.
\end{theorem}

\nid {\bf Proof.} Suppose $\psi:{\frak g}\to M$ is in the image of $\delta^0$ so that $\psi(g)=-gm$ for some fixed $m\in M$. The classical theory the guarantees us that $\delta_{\rm cl.}\psi=0$.  Moreover, for every $g\in{\frak g}$, we have
\[ \tilde\psi(g)=\psi(g^{[p]})-g^{p-1}\psi(g)=-g^{[p]}m+g^{p-1}gm=-(g^{[p]}-g^p)m=0 \] 
since $M$ is a restricted module. This shows that $\psi\in\Ker\delta^1$ and hence $H^1({\frak g};M)$ is well defined. If $\alpha\in H^1({\frak g};M)$ is a restricted cohomology class and $\psi\in\alpha$ represents $\alpha$, we must have $\delta^1\psi=(\delta_{\rm cl.}\psi,\tilde\psi)=(0,0)$ so that in particular, $\psi$ is a classical cocycle. Therefore $\psi$ also represents an ordinary cohomology class which we denote by $i(\alpha)\in H_{\rm cl.}^1({\frak g};M)$.  If $\psi^\prime\in \alpha$, then $\psi- \psi^\prime=\delta^0(m)=\delta^0_{\rm cl.}(m)$ for some $m\in M$ and hence $\psi$ and $\psi^\prime$ are cohomologous in the classical sense as well.  Therefore the mapping $i:H^1({\frak g};M)\to H_{\rm cl.}^1({\frak g};M)$ given by $i:\alpha\mapsto i(\alpha)$ is well defined and clearly linear. Finally, if $i(\alpha)=0$, then $\psi=\delta_{\rm cl.}^0(m)$ for some $m\in M$ and hence $\alpha=0$ in the restricted cohomology as well. Therefore the map $i$ is an injection and the proof of the theorem is complete. 
\qed

Theorem (\ref {thm:1}) gives an explicit description of $H^1(\frak g;M)$ as a subspace of $H_{\rm cl.}^1(\frak g;M)$.  Namely, a cohomology class $\alpha\in H_{\rm cl.}^1(\frak g;M)$ represents a restricted cohomology class if and only if $\tilde\psi=0$ for all representatives $\psi\in\alpha$.  We therefore can immediately investigate the restricted versions of the algebraic interpretations of $H_{\rm cl.}^1(\frak g;M)$.  We postpone these investigations until the subsequent section.

The situation for $C^3({\frak g};M)$ is considerably more complicated.  However, once we obtain the desired space and coboundary operator, we can make all necessary calculations of low dimensional cohomology and its algebraic interpretations.  That is, we are usually content with the beginning of the cochain complex
\[ 0\to C^0({\frak g};M)\to C^1({\frak g};M)\to C^2({\frak g};M)\to C^3({\frak g};M)
\] 
and we nearly have it. 

If $\alpha:\Lambda^3\frak g\to M$ is a skew-symmetric multilinear map on $\frak g$ and $\beta:\frak g\times\frak g\to M$ , we say that $\beta$ has the {\it $**$-property with respect to $\alpha$} if the following conditions hold:
\begin{itemize}
\item[(\romannumeral 1)] $\beta (g,h)$ is linear with respect to $g$.
\item[(\romannumeral 2)] $\beta (g,\lambda h)=\lambda^p\beta(g,h)$ for all $\lambda\in\mathbb F$.
\item[(\romannumeral 3)] \begin{eqnarray*}\label{2dstarprop}
\beta(g,h_1+h_2)& = & \beta(g,h_1)+\beta(g,h_2)- \\ 
 &  & \sum_{\stackrel{l_1,\dots,l_p=1 {\rm or} 2}{\scriptscriptstyle l_1=1, l_2=2}}\frac{1}{\#\{l_i=1\}}\sum_{j=0}^{p-2}(-1)^j\sum_{k=1}^j \bino{j}{k} h_{l_p}\cdots h_{l_{p-k-1}}\cdot \\
 & &\alpha ([g,h_{l_{p-k}},\cdots,h_{l_{p-j+1}}],[h_{l_1},\cdots,h_{l_{p-j-1}}],h_{l_{p-j}}). 
\end{eqnarray*}
\end{itemize}
Our space of 3-dimensional cochains is then defined as 
\[ C^3({\frak g};M)=\{(\alpha,\beta): \alpha\in C^3(\frak g;M), \beta:{\frak g}\times {\frak g}\to M\ \hbox {\rm has the $**$-property w.r.t. $\alpha$}\} \]
Again, it is evident that if $\beta$ and $\beta^\prime$ have the $**$-property with respect to $\alpha$ and $\alpha^\prime$ respectively, then $\beta+\beta^\prime$ has the $**$-property with respect to $\alpha+\alpha^\prime$, and hence $C^3(\frak g,M)$ is a vector space over $\mathbb F$ by pointwise addition in each coordinate.  As in the 2-dimensional case, given an element $\alpha\in C^3_{\rm cl.}(\frak g;M)$, there are numerous maps $\beta:\frak g\times\frak g\to M$ having the $**$-property with respect to $\alpha$.  Indeed, we may again assign the values of $\beta$ arbitrarily on a basis for $\frak g$ in each coordinate and conditions (\romannumeral 1)-(\romannumeral 3) above completely determine $\beta$ for a given $\alpha$.  We then have an exact sequence of $\mathbb F$-vector spaces
\[0\to K\to C^3(\frak g;M)\to C^3_{\rm cl.}(\frak g;M)\to 0\]
where the map $C^3(\frak g;M)\to C^3_{\rm cl.}(\frak g;M)$ is given by $(\alpha,\beta)\mapsto \alpha$ and $K$ denotes the kernel of this map.  This kernel consists of all pairs $(0,\beta)$ and $\beta$ has the $**$-property with respect to $0$ if and only if $\beta(g,h)$ is linear in $g$ and $p$-semilinear in $h$.  In this connection we denote the kernel $K$ by $\Hom_{\mathbb F}(\frak g\otimes \overline{\frak g};M)$.  If we again let $\dim_{\mathbb F}\frak g=n$ and $\dim_{\mathbb F}M=m$, then we see that 
\[\dim_{\mathbb F}C^3(\frak g;M)=n^2m+\bino{n}{3}m=\left (\frac{n(n+1)(n+2)}{6}\right )m\]
which agrees with our result in Corollary (\ref {cor:dimensions}).

Now, an element $(\phi,\omega)\in C^2(\frak g;M)$ induces a map $\beta:\frak g\times \frak g\to M$ by the formula
\begin{equation}
\label{3induced}
\beta(g,h)=\phi(g,h^{[p]})-\sum_{i+j=p-1}(-1)^ih^i\phi([g,\underbrace{h,\cdots,h}_j],h)+g\omega(h).
\end{equation}
To define the coboundary operator $\delta^2:C^2({\frak g};M)\to C^3({\frak g};M)$ we will need the following lemma.
\begin{lemma}
\label{starstarprop}
Given an element $(\phi,\omega)\in C^2(\frak g;M)$, the map $\beta$ defined in (\ref {3induced}) satisfies the $**$-property with respect to $\delta^2_{\rm cl.}\phi$.
\end{lemma}
The proof of Lemma (\ref {starstarprop}) is a computation in which we will use two combinatorial identities.  We therefore state and prove these results before giving the proof of Lemma (\ref {starstarprop}). 

\begin{lemma}
\label{combin1}
\[ \bino{p-1-s}{t}\equiv (-1)^{s+t}\bino{p-1-t}{s}\pmod {p}.\]
\end{lemma}

{\bf Proof.} The proof is a straightforward computation.  It is well known that $(p-1)!\equiv (-1) \pmod {p}$.  Therefore if $r\le p$, 
\[(p-r)!=\frac{(p-1)!}{(p-1)\cdots (p-r+1)}\equiv (-1)^{r-1}\frac{-1}{1\cdots (r-1)}=\frac{(-1)^r}{(r-1)!}\pmod {p}.\]
It follows that 
\begin{eqnarray*}
\lefteqn{\bino{p-1-s}{t}=\frac{(p-1-s)!}{t!(p-1-s-t)!}} \\
&\equiv & \frac{(-1)^{s+1}/s!}{((-1)^{t+1}/(p-1-t)!)(p-1-(s_t))!} \pmod {p}\\
&= & (-1)^{s+t}\frac{(p-1-t)!}{s!(p-1-s-t)!}=(-1)^{s+t}\bino{p-1-t}{s}.
\end{eqnarray*}
\qed

\begin{lemma}
\label{combin2}
If $a>b$, then 
\[ \sum_{i=b}^{a-c}(-1)^i\bino{a}{i+c}\bino{i}{b}=(-1)^{b}\bino{a-b-1}{c-1}.\]
\end{lemma} 

{\bf Proof.} From the identity $(1+t)^a(1+t)^{-b-1}=(1+t)^{a-b-1}$, we have 
\begin{equation}
\label{eq:coeff}
\sum_{i\ge 0}\bino{a}{i}t^i\cdot\sum_{j\ge 0}(-1)^j\bino{j+b}{b}t^j=\sum_{k\ge 0}\bino{a-b-1}{k}t^k.
\end{equation}
Equating the coefficients of $t^k$ in (\ref {eq:coeff}), we have
\begin{eqnarray*}
\bino{a-b-1}{k}&=& \sum_{i+j=k}(-1)^j\bino{a}{i}\bino{j+b}{b} \\
 & =& \sum_{i=0}^k(-1)^{k-i}\bino{a}{a-i}\bino{k+b-i}{b}\\
 & =& \sum_{i=b}^{k+b}(-1)^{i-b}\bino{a}{i+(a-b-k)}\bino{i}{b}
\end{eqnarray*}
The result follows immediately setting $a-b-k=c$. 
\qed

{\bf Proof of Lemma (\ref {starstarprop})} It is easy to see from the definition of $\beta$ that conditions (\romannumeral 1) and (\romannumeral 2) of the $**$-property are satisfied. It remains to verify condition (\romannumeral 3) with $\alpha=\delta^2_{\rm cl.}\phi$.  If $g,h_1,h_2\in\frak g$, we have 
\begin{eqnarray}
\lefteqn{\beta(g,h_1+h_2)=\underline {\phi(g,h_1^{[p]})}+\underline {\phi(g,h_2^{[p]})}} \nonumber \\
&+&\sum_{\stackrel{k_s=1\ {\rm or}\  2}{\scriptscriptstyle k_1=1,k_2=2}}\frac{1}{\#\{k_s=1\}}\phi(g,[h_{k_1},\dots,h_{k_p}]) \nonumber \\
&-& \sum_{i+j=p-1}\sum_{k_s=1\ {\rm or}\ 2} (-1)^i h_{k_1}\cdots h_{k_i}\phi([g,h_{k_{i+1}},\dots,h_{k_{p-1}}],h_{k_p})  \label{eq:star} \\
&+&  \underline {g\omega(h_1)} + \underline {g\omega(h_2)}  \nonumber \\
&+&\sum_{\stackrel{k_s=1\ {\rm or}\  2}{\scriptscriptstyle k_1=1,k_2=2}}\frac{1}{\#\{k_s=1\}}\sum_{j=0}^{p-2}(-1)^j gh_{k_p}\cdots h_{k_{p-j+1}}\phi([h_1,\dots,h_{k_{p-j-1}}],h_{k_{p-j}}).\nonumber 
\end{eqnarray}
The four underlined terms along with the summands from the term (\ref {eq:star}) with all $k_s=1$ or all $k_s=2$ together make up $\beta(g,h_1)+\beta(g,h_2)$.  Therefore we must verify the remaining terms account for the double sum in condition (\romannumeral 3) with $\alpha=\delta^2_{\rm cl.}\phi$.  Using the relation
\begin{eqnarray*}
\lefteqn{\delta_{\rm cl.}\phi(u,v,w)=}\\
& & \phi([uv],w)+\phi([vw],u)+\phi([wu],v)-u\phi(v,w)-v\phi(w,u)-w\phi(u,v),
\end{eqnarray*}
we rewrite each term
\[gh_{k_p}\cdots h_{k_{p-j+1}}\phi([h_1,\dots,h_{k_{p-j-1}}],h_{k_{p-j}})\]
as a sum of one term involving $\delta_{\rm cl.}\phi$ and five terms involving $\phi$.  The terms with $\delta_{\rm cl.}\phi$ make up the sum
\begin{equation}
-\sum_{\stackrel{k_s=1\ {\rm or}\  2}{\scriptscriptstyle k_1=1,k_2=2}}\frac{1}{\#\{k_s=1\}}\sum_{j=0}^{p-2}(-1)^j \sum_{\stackrel{A+B=}{\scriptscriptstyle\{p-j+1,\dots,p\}}} h_{k_A}\delta_{\rm cl.}\phi([gh_{k_B}],[h_{k_1},\dots,h_{k_{p-j-1}}],h_{k_{p-j}}),\nonumber 
\end{equation}
which is equal to
\begin{eqnarray*}
\lefteqn{-\sum_{\stackrel{l_1,\dots,l_p=1 {\rm or} 2}{\scriptscriptstyle l_1=1, l_2=2}}\frac{1}{\#\{l_i=1\}}\sum_{j=0}^{p-2}(-1)^j\sum_{k=1}^j \bino{j}{k} h_{l_p}\cdots h_{l_{p-k-1}}\cdot} \\
 &&\delta_{\rm cl.} ([g,h_{l_{p-k}},\cdots,h_{l_{p-j+1}}],[h_{l_1},\cdots,h_{l_{p-j-1}}],h_{l_{p-j}}). 
\end{eqnarray*} 
It remains to show that the rest of the right hand side of the expression for $\beta(g,h_1+h_2)$ cancels completely.  That is, we must show that the following sum is identically zero.
\begin{eqnarray}
\lefteqn{-\sum_{\stackrel{\stackrel{k_s=1\ {\rm or}\ 2}{\scriptscriptstyle{\rm not\  all}\  k_s=1}}{\scriptscriptstyle {\rm not\  all}\  k_s=2}}\sum_{\stackrel {i+j=p-1}{\scriptscriptstyle i>0}} (-1)^i h_{k_1}\cdots h_{k_i}\phi([g,h_{k_{i+1}},\dots,h_{k_{p-1}}],h_{k_p})} \label{eq:circ1} \\
&+& \sum_{\stackrel{k_s=1\ {\rm or}\  2}{\scriptscriptstyle k_1=1,k_2=2}}\frac{1}{\#\{k_s=1\}}\phi(g,[h_{k_1},\dots,h_{k_p}])\label{eq:circ2} \\
&+& \sum_{\stackrel{k_s=1\ {\rm or}\  2}{\scriptscriptstyle k_1=1,k_2=2}}\sum_{\stackrel {i+j=p-1}{\scriptscriptstyle i>0}}\sum_{\stackrel{A+B=}{\scriptscriptstyle\{p-j+1,\dots,p\}}}\frac{(-1)^j}{\#\{k_i=1\}}h_{k_A}\cdot \nonumber \\
 & & \bigg (\phi([[g,h_{k_B}],[h_{k_1},\dots,h_{k_i}]],h_{k_{i+1}}) \label{eq:circ3} \\
 &- & \phi([g,h_{k_B}],[h_{k_1},\dots,h_{k_{i+1}}])\label{eq:circ4} \\
 &- & \phi([g,h_{k_B},h_{k_{i+1}}],[h_{k_1},\dots,h_{k_i}])\label{eq:circ5} \\
 &+ & [h_{k_1},\dots,h_{k_i}]\phi([g,h_{k_B}],h_{k_{i+1}})\label{eq:circ6} \\
 &- & h_{k_{i+1}}\phi([g,h_{k_B}],[h_{k_1},\dots,h_{k_i}])\bigg ).\label{eq:circ7}
\end{eqnarray}
The entire sum consists of terms 
\[\underbrace{h_{i_1}\cdots h_{i_s}}_s\phi([g,\underbrace {h_{i_{s+1}},\dots,h_{i_{s+t}}}_t],[\underbrace {h_{i_{s+t+1}},\dots, h_{i_p}}_u])\]
where $s\ge 0$, $t\ge 0$, $u\ge 1$ and $s+t+u=p$.  We consider four cases.

{\bf Case 1.} $t\ge 1$ and $u\ge 2$.  These terms come only from (\ref {eq:circ4}),(\ref {eq:circ5}) and  (\ref {eq:circ7}).  By counting the number of occurrences in each line, we see the coefficients in (\ref {eq:circ4}),(\ref {eq:circ5}) and  (\ref {eq:circ7}) are respectively
\[-(-1)^{s+t}\frac{1}{\#(h_1)}\bino{s+t}{s}\]
\[-(-1)^{s+t-1}\frac{1}{\#(h_1)}\bino{s+t-1}{s}\]
\[-(-1)^{s+t-1}\frac{1}{\#(h_1)}\bino{s+t-1}{t}.\]
Therefore a well know identity from Pascal's triangle implies these terms cancel.

{\bf Case 2.} $s\ge 1$ and $t=0$.  These terms come only from (\ref {eq:circ4}) and  (\ref {eq:circ7}). The coefficients are respectively
\[-(-1)^{s}\frac{1}{\#(h_1)}\]
\[-(-1)^{s-1}\frac{1}{\#(h_1)}\]
so that the terms cancel.

{\bf Case 3.} $s=t=0$.  These terms come only from (\ref {eq:circ2}) and  (\ref {eq:circ4}). This time the coefficients are respectively
\[\frac{1}{\#(h_1)}\]
\[-\frac{1}{\#(h_1)}\]
so that, again,  the terms cancel.

{\bf Case 4.} $u=1$.  In this case, the terms have the form 
\begin{equation} 
\label{eq:dagger}
h_{i_1}\cdots h_{i_s}\phi([g,h_{i_{s+1}},\dots,h_{i_{p-1}}],h_{i_p}).
\end{equation} 
This term appears in (\ref {eq:circ1}),(\ref {eq:circ3}),(\ref {eq:circ5}),(\ref {eq:circ6}) and (\ref {eq:circ7}).  In (\ref {eq:circ1}), we note that not all the $i_t$ are the same and the coefficient is $-(-1)^s$.  We only have terms in (\ref {eq:circ5}) if $i=1$ and $i_p\ne 2$.  In this case $j=p-2$ so that $(-1)^j=-1$.  If $i_p=1$, then each pair $A,B$ with $|A|=s$ is counted once and the coefficient is 
\[\frac{1}{\#(h_1)}\bino{p-2}{s}=\frac{1}{\#(h_1)}(-1)^{s+1}\bino{p-s-1}{1}=\frac{1}{\#(h_1)}(-1)^s(s+1).\]
The only difference for (\ref {eq:circ7}) is that $|A|=s-1$ and hence the coefficient is
\[\frac{1}{\#(h_1)}\bino{p-2}{s-1}=\frac{1}{\#(h_1)}(-1)^{s-1}s.\]
Therefore adding the terms in (\ref {eq:circ5}) and (\ref {eq:circ7}) we have
\[(-1)^s\frac{1}{\#(h_1)}\cdot \left\{\begin{array}{cl}
1& \mbox {if}\  i_p=1 \\
0& \mbox {if}\  i_p=2.
\end{array}\right. \]
To investigate the coefficients in (\ref {eq:circ3}), we first note that Lemma (\ref {bracket}) implies that 
\begin{eqnarray*}
\lefteqn{(-1)^j h_{k_A}\phi([[g,h_{k_B}],[h_{k_1},\dots,h_{k_i}]],h_{k_{i+1}})=} \\
& & \sum_{C+D=\{2,\dots,i\}}(-1)^{j+|C|} h_{k_A}\phi([g,h_{k_B},h_{k_C},h_{k_1},h_{k_D}],h_{k_{i+1}}).
\end{eqnarray*}
Now, to make 
\[h_{k_A}\phi([g,h_{k_B},h_{k_C},h_{k_1},h_{k_D}],h_{k_{i+1}})\]
equal to 
\[h_{i_1}\cdots h_{i_s}\phi([g,h_{i_{s+1}},\dots,h_{i_{p-1}}],h_{i_p})\]
we must first choose $q$ ($s+1\le q\le p-1$) with $i_q=1$ to make $h_{i_q}=h_{k_1}$.  Then we should choose $j$ ($s\le j\le q$) to make $|A|+|B|=j$.  Then we should choose $A\subset \{p-j-1,\dots, p\}$ with $|A|=s$. Here $h_{k_A}$ should equal $h_{i_1}\cdots h_{i_s}$ ($i_1>\cdots > i_s$).  Finally we choose 
\[C\subset \{2,\dots,i-1=p-2-j\}\]
with $|C|=q-1-j$.  Here $h_{k_C}$ should be $h_{i_{j+1}},\dots,h_{i_{q-1}}$ ($i_{j+1}<\cdots i_{q-1}$).  Using Lemma (\ref {combin1}), the number of such choices is 
\[\bino{j}{s}\bino{p-2-j}{q-1-j}=\bino{j}{s}\bino{p-2-j}{p-1-q}\equiv (-1)^{q+j+1}\bino{j}{s}\bino{q}{j+1}\pmod {p}.\]
We remark here that 
\[(-1)^{j+|C|}=(-1)^{|A|+|B|+|C|}=(-1)^{q-1}.\]
Therefore the the coefficient in (\ref {eq:circ3}) is 
\begin{eqnarray*}
\lefteqn{\frac{1}{\#(h_1)}\sum_{\stackrel{s-1\le q\le p-1}{\scriptscriptstyle i_q=1}} (-1)^{q-1}\sum_{j=s}^{q-1} (-1)^{q+j+1}\bino{q}{j+1}\bino{j}{s}} \\
& & = \frac{1}{\#(h_1)} \sum_{\stackrel{s-1\le q\le p-1}{\scriptscriptstyle i_q=1}}\sum_{j=s}^{q-1}(-1)^j\bino{q}{j+1}\bino{j}{s} \\
& & =\frac{1}{\#(h_1)} \sum_{\stackrel{s-1\le q\le p-1}{\scriptscriptstyle i_q=1}}(-1)^s\bino{q-s-1}{0} \\
& & = (-1)^s \frac{1}{\#(h_1)}\#\{i_q=1 : s+1\le q\le p-1\}.
\end{eqnarray*}
To compute the coefficient in (\ref {eq:circ6}), we note that Lemma (\ref {bracket}) implies
\begin{eqnarray*}
\lefteqn{(-1)^j h_{k_A}[h_{k_1},\dots,h_{k_i}]\phi([g,h_{k_B}],h_{k_{i+1}})=} \\
& & \sum_{C+D=\{2,\dots,i\}}(-1)^{j+|C|} h_{k_A}h_{k_C}h_{k_1}h_{k_D}\phi([g,h_{k_B}],h_{k_{i+1}}).
\end{eqnarray*}
Now, to make 
\[h_{k_A}h_{k_C}h_{k_1}h_{k_D}\phi([g,h_{k_B}],h_{k_{i+1}})\]
equal to 
\[h_{i_1}\cdots h_{i_s}\phi([g,h_{i_{s+1}},\dots,h_{i_{p-1}}],h_{i_p})\]
we must first choose $q$ ($1\le q\le s$) with $i_q=1$ to make $h_{i_q}=h_{k_1}$.  Next we choose $a=|A|$ ($0\le a\le q-1$).  After this we make two choices: first choose $A$, with $|A|=a$, from a set of $|A|+|B|=a+(p-s-1)$ elements, and then choose $C$, with $|C|=q-1-a$, from a set of $|C|+|D|=s-a-1$ elements.  The total number of such choices is 
\begin{eqnarray*}
\bino{p-s-1+a}{a}\bino{s-1-a}{q-1-a}&=&\bino{p-1-s+a}{p-1-s}\bino{s-1-a}{s-q}\\
& \equiv& (-1)^a\bino{s}{s-a}\bino{s-a-1}{s-q}\pmod {p}.
\end{eqnarray*}
We remark again that 
\[(-1)^{j+|C|}=(-1)^{|A|+|B|+|C|}=(-1)^{q+s+1}\]
so that the total contribution of (\ref {eq:circ6}) modulo $p$ is 
\begin{eqnarray*}
\lefteqn{\frac{1}{\#(h_1)}\sum_{\stackrel{1\le q\le s}{\scriptscriptstyle i_q=1}} (-1)^{q+s+1}\sum_{a=0}^{q-1} (-1)^{a}\bino{s}{s-a}\bino{s-a-1}{s-q}} \\
& & = \frac{1}{\#(h_1)} \sum_{\stackrel{1\le q\le s}{\scriptscriptstyle i_q=1}}(-1)^q\sum_{t=s-q}^{s-1}(-1)^t\bino{s}{t+1}\bino{t}{s-q} \\
& & =\frac{1}{\#(h_1)} \sum_{\stackrel{1\le q\le s}{\scriptscriptstyle i_q=1}}(-1)^s\bino{q-1}{0} \\
& & = (-1)^s \frac{1}{\#(h_1)}\#\{i_q=1 : 1\le q\le s\}\\
& & =(-1)^s.
\end{eqnarray*}
If we add the coefficients from (\ref {eq:circ6}),(\ref {eq:circ3}),(\ref {eq:circ5}) and (\ref {eq:circ7}) we have 
\begin{eqnarray*}
\lefteqn{(-1)^s\frac{1}{\#(h_1)}\bigg (\#\{i_q=1:1\le q\le s\}+} \\
& & \#\{i_q=1:s+1\le q\le p-1\} + \#\{i_q=1: q=p\}\bigg ) \\
& & = (-1)^s\frac{\#(h_1)}{\#(h_1)}=(-1)^s.
\end{eqnarray*}
Therefore we see that these terms cancel (\ref {eq:circ1}) and the proof is complete.
\qed

Lemma (\ref {starstarprop}) allows to define the coboundary operator $\delta^2:C^2(\frak g;M)\to C^3(\frak g;M)$ by the formula
\[\delta^2:(\phi,\omega)\mapsto (\delta^2_{\rm cl.}\phi,\beta)\]
where $\beta$ is defined in (\ref {3induced}).  We then have the following theorem.
\begin{theorem}
\label{thm:2}
In the above notations, $\delta^2\delta^1=0$ so that the quotient $H^2(\frak g;M)=\Ker\delta^2/\Im\delta^1$ is well defined.
\end{theorem}

The proof of Theorem (\ref {thm:2}) is another computation which requires one more combinatorial lemma.

\begin{lemma}
\label{combin3}
If $2\le n\le p$, then 
\[\sum_{k=0}^{n-1}\bino{p-n+k}{k}=\bino{p}{n-1}\equiv 0\pmod {p}.\]
\end{lemma}

{\bf Proof.}  The equality $\bino{p}{n-1}\equiv 0\pmod {p}$ for $2\le n\le p$ is well known.  Using the identity
\[\bino{p}{k}=\bino{p-1}{k-1}+\bino{p-1}{k}\]
from Pascal's triangle, it is easy to show by induction on $n$ that
\begin{equation}
\label{eq:pascal}
\bino{p}{n-1}=\sum_{k=0}^{n-1}\bino{p-k-1}{n-k-1}.
\end{equation}
We remark that we define $\bino{s}{t}=0$ whenever $s<t$ for purposes of our induction. If we re-index in (\ref {eq:pascal}) with $\kappa=n-k-1$, then we have 
\begin{equation*}
\bino{p}{n-1}=\sum_{\kappa=0}^{n-1}\bino{p-n+\kappa}{\kappa}.
\end{equation*}
\qed

{\bf Proof of Theorem (\ref {thm:2}).} By definition, if $\psi\in C^1(\frak g;M)$, then $\delta^2\delta^1(\psi)=(\delta^2_{\rm cl.}\delta^1_{\rm cl.}(\psi),\beta)$ where $\beta:\frak g\times \frak g:\to M$ is defined by
\begin{equation}
\label{eq:pound}
\beta(g,h)=\delta^1_{\rm cl.}\psi(g,h^{[p]})-\sum_{i+j=p-1}(-1)^ih^i\delta^1_{\rm cl.}\psi([g,\underbrace{h,\cdots,h}_j],h)+g\psi(h^{[p]})-gh^{p-1}\psi(h).
\end{equation}
Of course $\delta^2_{\rm cl.}\delta^1_{\rm cl.}(\psi)=0$ so that it remains to show that $\beta=0$.  Using the identity 
\begin{equation}
\label{eq:swap}
\delta^1_{\rm cl.}\psi(g,h)=-g\psi(h)+h\psi(g)+\psi([gh]),
\end{equation}
we rewrite each term in (\ref {eq:pound}) involving $\delta^1_{\rm cl.}\psi$ as a sum of three terms involving $\psi$.  Expanding we have
\begin{eqnarray}
\lefteqn{\beta(g,h)=\underline{-g\psi(h^{[p]})}+\underline {\underline{h^{[p]}\psi(g)}}+\underline {\underline {\psi([gh^{[p]}])}}}\nonumber \\
& & +\sum_{i+j=p-1}(-1)^ih^i[g,\underbrace{h,\cdots,h}_j]\psi(h)\label{eq:A} \\
& & -\sum_{i+j=p-1}(-1)^ih^{i+1}\psi([g,\underbrace{h,\cdots,h}_j])\label{eq:B} \\
& & -\sum_{i+j=p-1}(-1)^ih^{i}\psi([g,\underbrace{h,\cdots,h}_{j+1}])\label{eq:C} \\
& & +\underline{g\psi(h^{[p]})}-gh^{p-1}\psi(h).\label{eq:D}
\end{eqnarray}
Clearly the two underlined terms cancel.  Moreover, if $i<p-1$, then every term in (\ref {eq:B}) appears in (\ref {eq:C}) with the opposite sign.  This leaves on the term $j=p-1$ in (\ref {eq:C}).  If $i=p-1$, the remaining term in (\ref {eq:B}) is $-h^p\psi(g)$ ($[g]=g$), and if $j=p-1$, the remaining term in (\ref {eq:C}) is $-\psi([g,\underbrace{h,\cdots,h}_{p}]))$.  Therefore these terms cancel the double underlined terms since $h^{[p]}\psi(g)=h^p\psi(g)$ and $[g,\underbrace{h,\cdots,h}_{p}]=[gh^{[p]}]$.  Now, recalling that 
\[[g,\overbrace {h,\dots,h}^j]=\sum_{s+t=j}(-1)^{s}\bino{j}{s}h^sgh^t, \]
we expand the bracket in (\ref {eq:A}) and we have
\begin{eqnarray*}
\lefteqn{\beta(g,h)=\sum_{i+j=p-1}\sum_{s+t=j}(-1)^{i+s}\bino{j}{s} h^{i+s} g h^t\psi(h)}\\
& & -gh^{p-1}\psi(h) \\
& & =gh^{p-1}\psi(h) +\sum_{\stackrel {s+t=p-1}{\scriptscriptstyle s>0}}(-1)^s\bino{p-1}{s} h^s g h^t \psi(h)\\
& & +\sum_{\stackrel{i+j=p-1}{\scriptscriptstyle i>0}}\sum_{s+t=j}(-1)^{i+s}\bino{j}{s} h^{i+s} g h^t\psi(h)\\
& & -gh^{p-1}\psi(h).
\end{eqnarray*}
If we rearrange the remaining terms, we have 
\[\beta(g,h)=\sum_{i=2}^p(-1)^{i-1}\left [\sum_{j=0}^{i-1} \bino{p-i+j}{j}\right ]h^{i-1}gh^{p-i}\psi(h).\]
But by Lemma (\ref {combin3}), each coefficient in this sum is zero modulo $p$ so that $\beta=0$ as claimed.
\qed

We remark here that the map $C^2(\frak g;M)\to C^2_{\rm cl.}(\frak g;M)$ given by $(\phi,\omega)\mapsto \phi$ is obviously a cochain map so that we have  a map $H^2(\frak g;M)\to H^2_{\rm cl.}(\frak g;M)$ that sends the cohomology class of $(\phi,\omega)$ to the (classical) cohomology class of $\phi$.  Indeed, we have the following commutative diagram
\begin{equation}
\label{ladder}
\begin{array}{ccccccccc}
0&\longrightarrow&C^0(\frak g;M)&\longrightarrow&C^1(\frak g;M) &\longrightarrow&C^2(\frak g;M) & \longrightarrow & C^3(\frak g;M) \\
& & \|& \circ & \| & \circ & \downarrow &\circ & \downarrow \\
0&\longrightarrow&C_{\rm cl.}^0(\frak g;M)&\longrightarrow&C_{\rm cl.}^1(\frak g;M) &\longrightarrow&C_{\rm cl.}^2(\frak g;M) & \longrightarrow & C_{\rm cl.}^3(\frak g;M)
\end{array}
\end{equation}
where the horizontal maps are the coboundary operators.  This diagram shows that we have a map $H^k(\frak g;M)\to H_{\rm cl.}^k(\frak g;M)$ for $k\le 2$.  Theorem (\ref {thm:1}) states that this map is injective for $k=1$; this is obvious from the left most square in (\ref {ladder}).  The map $H^2(\frak g;M)\to H_{\rm cl.}^2(\frak g;M)$ fails to be injective.  Indeed, if $(\phi,\omega)\in C^2(\frak g;M)$ represents a restricted cohomology class and $\phi\in C^2_{\rm cl.}(\frak g;M)$ is cohomologous to zero; i.e. $\phi=\delta^1_{\rm cl.}\psi$,  it need not follow that $\omega=\tilde\psi$ so that $(\phi,\omega)$ is not necessarily cohomologous to zero as a restricted cocycle.


\section{Algebraic Interpretations}
In this section, we develop the analogs of the algebraic interpretations of low dimensional cohomology described in section (2.1.4) for restricted Lie algebras and show that equivalence classes of these objects are naturally encoded in the (restricted) cohomology spaces defined above.  In each example, we use the interpretations in the classical case along with our description of restricted cohomology to define corresponding notions for restricted Lie algebras.  In the case of 1-dimensional cohomology, we employ the injection $H^1(\frak g;M)\to H^1_{\rm cl.}(\frak g;M)$ of Theorem (\ref {thm:1}) to motivate appropriate definitions of restricted derivations and extensions of restricted modules.  In dimension zero, of course we still have $H^0(\frak g;M)=M^\frak g$ since this is the case for $H^0_{\rm cl.}(\frak g;M)$.  

We begin with the notion of a restricted derivation of a restricted Lie algebra $\frak g$.  Recalling that $H^1_{\rm cl.}(\frak g;\frak g)$ equal to the space $\Der(\frak g)/\ad (\frak g)$ of outer derivations of $\frak g$, Theorem (\ref {thm:1}) motivates the following definition.
\begin{definition}
\label{def:resderivation}
If $\frak g$ is a restricted Lie algebra over $\mathbb F$, an $\mathbb F$-linear map $D:\frak g\to \frak g$ is called a {\it restricted derivation of}\ $\frak g$ if for all $g,h\in\frak g$, 
\begin{itemize}
\item[{\rm (\romannumeral 1)}] $D([gh])=[gD(h)]+[D(g)h]$.
\item[{\rm (\romannumeral 2)}] $D(g^{[p]})=(\ad g)^{p-1}D(g)$.
\end{itemize}
\end{definition}
We denote the set of all restricted derivations of $\frak g$ by $\Der_{\rm res.}(\frak g)$; it is an $\mathbb F$-submodule of $\frak {gl}(\frak g)=\Hom_{\mathbb F}(\frak g,\frak g)$.  
\begin{lemma}
\label{thm:adjoint}
For each $g\in \frak g$, the map $\ad g:\frak g\to \frak g$ defined by
\[ \ad g :h\mapsto [gh]\]
is a restricted derivation of $\frak g$ and the map $\ad:\frak g\to \Der_{\rm res.}(\frak g)$ that sends $g$ to $\ad g$ is a restricted Lie algebra homomorphism. 
\end{lemma}

{\bf Proof.} It is well known that $\ad g$ is an ordinary derivation on $\frak g$ in the sense of (\romannumeral 1) of Definition (\ref {def:resderivation}).  If $h\in\frak g$, we have 
\[\ad g(h^{[p]})=[gh^{[p]}]=[g,h,\underbrace{h,\dots,h}_{p-1}]=(-1)^{p-1}(\ad h)^{p-1}(\ad g (h)).\]
But $(-1)^{p-1}=1$ so that $\ad g$ is a restricted derivation of $\frak g$.  Of course $\ad:\frak g\to\frak {gl}(\frak g)$ is an ordinary Lie algebra homomorphism, and since 
\[\ad g^{[p]}(h)=-[hg^{[p]}]=-[h,\underbrace{g,\dots,g}_{p}]=(-1)^{p+1}(\ad g)^p(h),\]
we see that $\ad$ is in fact a restricted Lie algebra homomorphism. 
\qed

Lemma (\ref {thm:adjoint}) gives $\frak g$ the structure of a $\frak g$-module and the notation $C^k(\frak g;\frak g)$ will always mean this particular module structure.  Borrowing terminology from the classical case, we will call restricted derivations of the form $\ad g$ {\it inner}\  and elements of the $\mathbb F$-module quotient $\Der_{\rm res.}(\frak g)/\ad(\frak g)$ {\it outer}.  We remark that actually, $\Der_{\rm res.}(\frak g)$ is a restricted Lie subalgebra of $\frak {gl}(\frak g)$ and that the equation
\[ [D,\ad g]=\ad D(g)\]
shows that $\ad(\frak g)$ is an ideal in $\Der_{\rm res.}(\frak g)$. Therefore we can view the quotient $\Der_{\rm res.}(\frak g)/\ad(\frak g)$ as a restricted Lie algebra.  With these definitions, we have the following theorem immediately.
\begin{theorem}
\label{thm:outerderivations}
The space $\Der_{\rm res.}(\frak g)/\ad(\frak g)$ of restricted outer derivations of $\frak g$ is equal to $H^1(\frak g;\frak g)$.
\end{theorem}

{\bf Proof.} If $\psi\in C^1(\frak g;\frak g)$, then $\delta^1\psi=(\delta^1_{\rm cl.}\psi,\tilde \psi)$ is a cocycle if and only if $\delta^1_{\rm cl.}\psi=0$ and $\tilde\psi=0$. Easily $\delta^1_{\rm cl.}\psi=0$ if and only if $\psi$ satisfies (\romannumeral 1) of Definition (\ref {def:resderivation}) and $\tilde\psi=0$ if and only if $\psi$ satisfies (\romannumeral 2) so that $Z^1(\frak g;\frak g)=\Der_{\rm res.}(\frak g)$.  Moreover, since $\delta^0=\delta^0_{\rm cl.}$, it follows that $\Im\delta^0=\ad(\frak g)$ and the proof is complete.
\qed

We note that Theorem (\ref {thm:outerderivations}) together with the remarks after Lemma (\ref {thm:adjoint}) imply that $H^1(\frak g;\frak g)$ is a restricted Lie algebra.  As in the case of classical Lie algebra cohomology, the space $H^1(\frak g;\frak g)$ has another interpretation.  We begin with a slightly more general notion of which this is a special case.  
\begin{definition} 
If $M$ and $N$ are restricted $\frak g$-modules, then a restricted extension of $N$ by $M$ is an exact sequence 
\begin{equation}
\label {extension}
0\to M\stackrel{\iota}{\longrightarrow} E \stackrel{\pi}{\longrightarrow} N\to 0
\end{equation}
of restricted $\frak g$-modules and homomorphisms.
\end{definition}  
Two extensions of $N$ by $M$ are equivalent if they can be included in a commutative diagram
\begin{equation} 
\label{equivextensions}
\begin{array}{ccccccccc}
 0 & \longrightarrow & M & \stackrel{\iota_1}{\longrightarrow}& E_1& \stackrel{\pi_1}{\longrightarrow} &N& \longrightarrow &0 \\
&&\|&\circ&\downarrow&\circ&\|& \\
0&\longrightarrow& M&\stackrel{\iota_2}{\longrightarrow}&E_2&\stackrel{\pi_2}{\longrightarrow}&N&\longrightarrow&0 
\end{array}
\end{equation}
We denote the set of equivalence classes of extensions of $N$ by $M$ by $\Ext(N,M)$.  We will need the following lemma.
\begin{lemma}
\label{hommodule}
If $M$ and $N$ are restricted $\frak g$-modules, then the space of $\mathbb F$-linear maps $\Hom_{\mathbb F}(N,M)$ is a restricted $\frak g$-module where we define for $g\in\frak g$, $\phi\in\Hom_{\mathbb F}(N,M)$ and $n\in N$
\[(g\phi)(n)=g\phi(n)-\phi(gn).\]
\end{lemma}

{\bf Proof.} It is well known that this action gives $\Hom_{\mathbb F}(N,M)$ the structure of an ordinary $\frak g$-module. Therefore it remains to verify that $g^{[p]}=g^p$ as operators on $\Hom_{\mathbb F}(N,M)$.  For a fixed $g\in \frak g$, we define two endomorphisms, $u_g$ and $v_g$, of $\Hom_{\mathbb F}(N,M)$
by $u_g(\phi)=g\circ\phi$ and $v_g=\phi\circ g$ where we consider $g\in\frak g$ as an operator on $M$ and $N$ respectively.  In this notation, the $\frak g$-module structure on $\Hom_{\mathbb F}(N,M)$ is given by $g\mapsto u_g-v_g$.  Now, as endomorphisms of $\Hom_{\mathbb F}(N,M)$, one can easily verify that $u_g$ and $v_g$ commute so that $(u_g-v_g)^p=u_g^p-v_g^p$.  Now, if $g\in\frak g$, $\phi\in \Hom_{\mathbb F}(N,M)$ and $n\in N$, we compute
\begin{eqnarray*}
(g^{[p]}\phi)(n)& = & g^{[p]}\phi(n)-\phi(g^{[p]}n) \\
  & = & g^p\phi(n)-\phi(g^pn) \\
 & = & [(u_g^p-v_g^p)(\phi)](n) \\
 & = & [(u_g-v_g)^p(\phi)](n) \\
 & = & [g^p\phi](n)
\end{eqnarray*}
and hence $\Hom_{\mathbb F}(N,M)$ is a restricted $\frak g$-module as claimed.
\qed

Given an extension (\ref {extension}) of $N$ by $M$, we construct an element of $C^1(\frak g;\Hom_{\mathbb F}(N,M))$ as follows.  Choose a $\mathbb F$-linear map $\rho:N\to E$ such that $\pi\rho=1_N$. We then define an element $\psi\in C^1(\frak g;\Hom_{\mathbb F}(N,M))$ by the formula
\[\psi(g)(n)=g\rho(n)-\rho(gn).\]
We note that, in fact, $g\rho(n)-\rho(gn)\in M$ since 
\[\pi(g\rho(n)-\rho(gn))=g\pi\rho(n)-\pi\rho(gn)=gn-gn=0.\]
We claim that the map $\psi$ is a restricted cocycle whose cohomology class depends only on the equivalence class of the extension (\ref {extension}).  Indeed, it follows from the classical cohomology theory that $\delta^1_{\rm cl.}\psi=0$ and that the classical cohomology class of $\psi$ depends only on the equivalence class of the extension (\ref {extension}).  Since $\delta^0=\delta^0_{\rm cl.}$, it remains only to verify that $\tilde\psi=0$.  Now, $\rho\in\Hom_{\mathbb F}(N,E)$ which is a restricted $\frak g$-module and, by definition, we have $\psi(g)=g\rho$ as maps $N\to M$.  Therefore 
\[g^{p-1}\psi(g)=g^p\rho=g^{[p]}\rho=\psi(g^{[p]})\]
so that $\tilde\psi(g)=\psi(g^{[p]})-g^{p-1}\psi(g)=0$ and our claim is established.  Conversely, if $\psi\in C^1(\frak g;\Hom_{\mathbb F}(N,M))$ is a restricted cocycle, we construct an extension of $N$ by $M$ as follows.  As a vector space, we set $E=N\oplus M$ with $\iota$ and $\pi$ the canonical inclusion and projection respectively.  For each $g\in \frak g$, we define an endomorphism of $E$ by 
\[g(n,m)=(gn,gm+\psi(g)(n)).\]
We claim that this gives $E$ the structure of a restricted $\frak g$-module and that the sequence $0\to M\to E\to N\to 0$ is exact in the category of restricted $\frak g$-modules.  The linearity conditions clearly hold so that we need only check the restricted module conditions for the pairing $(g,(n,m))$.  In fact, the same construction, {\it mutatis mutandis}, works for ordinary Lie algebra modules and the condition $[gh]=gh-hg$ is verified there.  We therefore will only check that $g^{[p]}=g^p$ as operators on $E$.  If $g\in\frak g$, then using the cocycle condition on $\psi$ along with the fact that $N$ and $M$ are restricted $\frak g$-modules, we have 
\[g^{[p]}(n,m)=(g^{[p]}n,g^{[p]}m+\psi(g^{[p]})(n))=(g^pn,g^pm+(g^{p-1}\psi(g))(n)).\]
On the other hand, a direct computation shows that 
\[g^p(n,m)=(g^pn,g^pm+\sum_{i+j=p-1} g^j (\psi(g)(g^in))).\]
Using the notation of Lemma (\ref{hommodule}), we recall that the commuting endomorphisms $u_g$ and $v_g$ of $\Hom_{\mathbb F}(N,M)$ satisfy the identity
\[(u_g-v_g)^{p-1}=\sum_{i+j=p-1} u_g^jv_g^i.\]
Moreover, we have $ (u_g^jv_g^i\psi)(n)=g^j\psi(g^in)$ for all $n\in N$ so that \begin{eqnarray*}
(g^{p-1}\psi(g))(n)& =&((u_g-v_g)^{p-1}\psi(g))(n)\\
 & =& \left (\sum_{i+j=p-1}u_g^jv_g^i\psi(g)\right )(n)\\
 & =& \sum_{i+j=p-1} g^j (\psi(g)(g^in)).
\end{eqnarray*}
Therefore $g^{[p]}=g^p$ as operators on $E$, and hence $E$ is a restricted module.  Evidently the canonical inclusion and projection are restricted Lie algebra homomorphisms with this module structure so that we have an extension of $N$ by $M$.  We claim that the equivalence class of the extension constructed above depends only on the cohomology class of $\psi$ so that we have an assignment $H^1(\frak g;\Hom_{\mathbb F}(N,M))\to \Ext(N,M)$.  Suppose that $\psi_1$ and $\psi_2$ are cohomologous and let $E_1$ and $E_2$ denote the corresponding extensions of $N$ by $M$.  If $f:N\to M$ satisfies $\delta^0 f=\psi_1-\psi_2$, then we define a map $E_1\to E_2$ by
\[
(n,m)\mapsto (n,m-f(n)).
\]
Clearly this map is an isomorphism of vector spaces making the diagram (\ref {equivextensions}) commute.  Therefore we need only verify that it commutes with the action of $\frak g$.  We  have
\[g(n,m)=(gn,gm+\psi_1(g)(n))\mapsto (gn,gm+\psi_1(g)(n)-f(gn))\]
and 
\[g(n,m-f(n))=(gn,gm-gf(n)+\psi_2(g)(n)).\]
Taking the difference of the second factors we have
\begin{eqnarray*}
\lefteqn{gm+\psi_1(g)(n)-f(gn)-gm+gf(n)-\psi_2(g)(n)}\\
&  &= (\psi_1(g)-\psi_2(g))(n)-(-gf(n)+f(gn))\\
& & = (\psi_1(g)-\psi_2(g))(n)-(\delta^0 f(g))(n)\\
& & =0
\end{eqnarray*}
Therefore the extensions $E_1$ and $E_2$ are equivalent and our assignment is well defined.  It is obvious that our constructions are inverse to each other so that we have shown the following theorem.
\begin{theorem}
\label{thm:rext}
If $M$ and $N$ are restricted $\frak g$-modules, then the set $\Ext(N,M)$ of equivalence classes of restricted extensions of $N$ by $M$ is in one to one correspondence with $H^1(\frak g;\Hom_{\mathbb F}(N,M))$.
\qed
\end{theorem} 
In particular, if $N=\mathbb F$ is regarded as a trivial $\frak g$-module, we recover the restricted version of one dimensional right extensions of the module $M$. 

\begin{corollary}
$H^1(\frak g;M)$ is in one to one correspondence with equivalence classes of one dimensional right extensions of the restricted $\frak g$-module $M$.
\qed
\end{corollary}

\begin{corollary}
$H^1(\frak g,\frak g)$ is in one to one correspondence with equivalence classes of one dimensional right extensions of the restricted Lie algebra $\frak g$.
\qed
\end{corollary}
To interpret these results in terms of Theorem (\ref {thm:1}), we note that the ordinary cocycle condition is sufficient to make $E$ an ordinary $\frak g$-module; whereas the restricted cocycle condition $\psi(g^{[p]})=g^{p-1}\psi(g)$ is necessary for $E$ to be a {\it restricted}\ $\frak g$-module.  This is evident in the proof given above.  If let let $\Ext_{\rm cl.}(N,M)$ denote the set of equivalence classes of extensions of $N$ by $M$ as ordinary $\frak g$-modules, then Theorem (\ref {thm:rext}) can be stated as follows.
\begin{theorem}
The canonical isomorphism $\Ext_{\rm cl.}(N,M)\stackrel{\sim}{\longrightarrow} H^1_{\rm cl.}(\frak g;\Hom_{\mathbb F}(N,M))$ maps $\Ext(N,M)$ onto $H^1(\frak g;\Hom_{\mathbb F}(N,M))$.
\qed
\end{theorem}
We now turn our attention towards algebraic interpretations of $H^2(\frak g;M)$.  Unlike the one dimensional case, the canonical map $H^2(\frak g;M)\to H^2_{\rm cl.}(\frak g;M)$ is not injective so that we cannot simply investigate a particular subspace of $H^2_{\rm cl.}(\frak g;M)$.  We begin with the notion of restricted extensions.  Following Hochschild in \cite {H}, we say that a restricted Lie algebra $\frak h$ is {\it strongly abelian}\ if in addition to $[\frak h\frak h]=0$, we also have $\frak h^{[p]}=0$.  We then make the following definition.
\begin{definition}
If $\frak g$ is a restricted Lie algebra and $\frak h$ is a strongly abelian restricted Lie algebra, then a restricted  extension of $\frak g$ by $\frak h$ is an exact sequence 
\begin{equation}
\label{rescentralext}
0\longrightarrow\frak h\stackrel{\iota}{\longrightarrow}\frak e\stackrel{\pi}{\longrightarrow}\frak g\longrightarrow 0
\end{equation}
of restricted Lie algebras and their homomorphisms. 
\end{definition} 
Two restricted extensions of $\frak g$ by $\frak h$ are equivalent if they can be included in the usual commutative diagram.  That is $\frak e_1$ is equivalent to $\frak e_2$ if there is a restricted Lie algebra isomorphism $\alpha:\frak e_1\to\frak e_2$ that fixes $\frak h$ elementwise and $\pi_2\alpha=\pi_1$.  We note that a restricted extension of $\frak g$ by $\frak h$ gives $\frak h$ the structure of a $\frak g$-module by the action 
\[g\cdot h=[\tilde gh]\]
where $\tilde g\in\frak e$ is any element satisfying $\pi(\tilde g)=g$.   This action is well defined because $\frak h$ is abelian.  Moreover, since $\pi(\tilde g^{[p]})=\pi(\tilde g)^{[p]}=g^{[p]}$, it is easy to see that $\frak h$ is a restricted module.  We remark here that if $\frak h$ is contained in the center of $\frak e$, then $\frak h$ is a trivial $\frak g$-module.  Such an extension is called central.  Given an extension (\ref {rescentralext}) of $\frak g$ by $\frak h$, we construct an element of $H^2(\frak g;\frak h)$ as follows.  Choose an element $\sigma\in\Hom_{\mathbb F}(\frak g,\frak e)$ such that $\pi\sigma=1_{\frak g}$ and define $\phi:\frak g\times \frak g\to\frak h$ and $\omega:\frak g\to \frak h$ by the formulae
\begin{eqnarray*}
\phi(g,g')&=&[\sigma(g)\sigma(g')]-\sigma[gg']\\
\omega(g) &=& \sigma(g)^{[p]}-\sigma(g^{[p]}).
\end{eqnarray*}
We note that since $\pi\sigma=1_{\frak g}$, we have $\Im\phi\subset\frak h$ and $\Im\omega\subset\frak h$.  Moreover, $\phi$ is clearly $\mathbb F$-bilinear and skew-symmetric so that $\phi\in C^2_{\rm cl.}(\frak g;\frak h)$.  We need the following lemma.
\begin{lemma}
In the above notations, the map $\omega$ has the $*$-property with respect to $\phi$ so that the pair $(\phi,\omega)$ is an element of $C^2(\frak g;\frak h)$.  Moreover $\delta^2(\phi,\omega)=0$ so that $(\phi,\omega)$ represents a restricted cohomology class in $H^2(\frak g;\frak h)$.
\end{lemma}

{\bf Proof.}  Obviously we have $\omega(\lambda g)=\lambda^p\omega(g)$ for all $g\in \frak g$ and all $\lambda\in\mathbb F$.  Moreover, it is evident from the definition of $\omega$ that 
\begin{eqnarray}
\lefteqn{\omega(g+g')=\omega(g)+\omega(g')}\nonumber \\
 & & +\sum_{\stackrel{g_j=g\ {\rm or}\  g'}{\scriptscriptstyle g_1=g,g_2=g'}} \frac{1}{\#(g)}\left ([\sigma(g_1),\dots,\sigma(g_p)]-\sigma([g_1,\dots,g_p])\right ).\label{pound}
\end{eqnarray}
Using the formula
\[\sigma[gg']=[\sigma(g)\sigma(g')]-\phi(g,g'),\]
and the definition of the $\frak g$-module structure on $\frak h$, we expand the term $\sigma([g_1,\dots,g_p])$.  We have
\begin{eqnarray*}
\lefteqn{\sigma([g_1,\dots,g_p])=[\sigma(g_1),\dots,\sigma(g_p)]}\\
& & +\sum_{k=0}^{p-2}(-1)^{k+1}g_p\cdots g_{p-k+1}\phi([g_1,\cdots,g_{p-k-1}],g_{p-k}).
\end{eqnarray*}
Substituting this last expression into (\ref {pound}) shows that $\omega$ has the $*$-property with respect to $\phi$.  We remark that our construction of $\phi$ is valid for ordinary extensions of $\frak g$ by $\frak h$ and it is well known that $\delta^2_{\rm cl.}\phi=0$.  Therefore the proof of the lemma is complete upon showing that the induced map $\beta$ in (\ref {3induced}) is identically zero.  Using the definitions of $\phi$, $\omega$ and the $\frak g$-module structure on $\frak h$, we have for all $g,g'\in\frak g$, 
\begin{eqnarray*}
\beta(g,g')&=&[\sigma(g)\sigma({g'}^{[p]})]-\sigma([g{g'}^{[p]}])\\
&-&\sum_{i+j=p-1}(-1)^i\underbrace{[\sigma(g'),\dots,\sigma(g')}_i,\phi([g,\underbrace{g',\cdots,g'}_j],g')]\\
&+&[\sigma(g)\sigma(g')^{[p]}]-[\sigma(g)\sigma({g'}^{[p]})].
\end{eqnarray*}
Now, $[g{g'}^{[p]}]=[g,\underbrace{g',\dots,g'}_p]$ so that using the identity
\[\sigma[gg']=[\sigma(g)\sigma(g')]-\phi(g,g'),\]
we expand $\sigma([g{g'}^{[p]}])$ and see it cancels each term in the middle sum leaving only the term 
\[ [\sigma(g),\underbrace{\sigma(g'),\dots,\sigma(g')}_p]=[\sigma(g)\sigma(g')^{[p]}].\]
Therefore we see that $\beta=0$ and the proof is complete.
\qed

Suppose that $\sigma':\frak g\to\frak e$ is another $\mathbb F$-linear splitting map and let $(\phi',\omega')$ denote the corresponding 2-cocycle constructed above.  If we let $\psi=\sigma-\sigma'$, then easily $\Im\psi\subset\frak h$ so that $\psi\in C^1(\frak g;\frak h)$.  Moreover it is well known that $\delta^1_{\rm cl.}\psi=\phi'-\phi$.  We claim that $\tilde\psi=\omega'-\omega$ so that $(\phi,\omega)$ and $(\phi',\omega')$ are cohomologous as restricted 2-cocycles.  To verify our claim, we first compute for any $g\in\frak g$:
\begin{eqnarray}
\tilde\psi(g)& =& \psi(g^{[p]})-g^{p-1}\psi(g) \nonumber \\
 & =& \psi(g^{[p]})-[\underbrace{\sigma(g)[\cdots [\sigma(g)}_{p-1}\psi(g)]\cdots ]]\nonumber \\
& =& \sigma(g^{[p]})-\sigma'(g^{[p]})-[\underbrace{\sigma(g)[\cdots [\sigma(g)}_{p-1}(\sigma(g)-\sigma'(g))]\cdots ]]\nonumber \\
& =& \sigma(g^{[p]})-\sigma'(g^{[p]})+[\underbrace{\sigma(g)[\cdots [\sigma(g)}_{p-1}\sigma'(g)]\cdots ]]\label{eq:lastline}
\end{eqnarray}
On the other hand we have 
\begin{eqnarray}
(\omega'-\omega)(g)& =& \sigma'(g)^{[p]}-\sigma'(g^{[p]})-\sigma(g)^{[p]}+\sigma(g^{[p]})\nonumber \\
& =&\sigma(g^{[p]})-\sigma'(g^{[p]})+\sigma'(g)^{[p]}-\sigma(g)^{[p]}\label{eq:lastlastline}
\end{eqnarray}
Now, since $\frak h$ is strongly abelian, $\psi(g)^{[p]}=0$ so that $(\sigma(g)-\sigma'(g))^{[p]}=0$.  It follows that 
\[\sigma'(g)^{[p]}-\sigma(g)^{[p]}=\sum_{\stackrel{g_j=\sigma(g)\ {\rm or}\  \sigma'(g)}{\scriptscriptstyle g_1=\sigma(g),g_2=\sigma'(g)}} \frac{1}{\#(\sigma(g))} [g_1,g_2,\cdots,g_p].\]
But $\sigma(g)$ and $\sigma'(g)$ are equal as operators on $\frak h$ so that for a fixed number of $\sigma(g)$, this term occurs $\#(\sigma(g))$ times.  This together with an inspection of equations (\ref {eq:lastline}) and (\ref {eq:lastlastline}) establishes our claim. 

Now, if $\frak e_1$ and $\frak e_2$ are equivalent restricted extensions, and $\alpha:\frak e_1\to\frak e_2$ is the map realizing this equivalence, then we choose two splitting maps $\sigma_i:\frak g\to\frak e_i$ and construct corresponding cohomology classes containing $(\phi_1,\omega_1)$ and $(\phi_2,\omega_2)$.  The same arguments given above show that the map $\psi=\alpha\sigma_1-\sigma_2$ takes values in $\frak h$ and $\delta^1\psi=(\phi_2-\phi_1,\omega_2-\omega_1)$.  Therefore we have a well defined map from the set of equivalence classes of restricted extensions of $\frak g$ by $\frak h$ to $H^2(\frak g;\frak h)$.  Conversely, if $(\phi,\omega)\in C^2(\frak g;\frak h)$ is a cocycle, we construct a restricted extension of $\frak g$ by $\frak h$ as follows.  We make $\frak h$ into a strongly abelian restricted Lie algebra by declaring that $[hh']=0$ and $h^{[p]}=0$ for all $h,h'\in\frak h$. We define $\frak e=\frak h\oplus \frak g$ as a vector space and we define the Lie bracket and $p$-operator in $\frak e$ by the formulae
\begin{eqnarray}
[(h,g)(h',g')]&=&(\phi(g,g'),[gg'])\label{jacobiext}\\
(h,g)^{[p]}&=&(\omega(g),g^{[p]}).\label{pext}
\end{eqnarray}
The bracket (\ref {jacobiext}) is clearly bilinear and skew-symmetric and it is well known that the Jacobi identity for (\ref {jacobiext}) is equivalent to $\delta^2_{\rm cl.}\phi=0$.  Moreover, the operation (\ref {pext}) is a $p$-operator precisely because $\omega$ has the $*$-property with respect to $\phi$ and the induced map $\beta$ in (\ref {3induced}) is identically zero.  Finally if $(\phi',\omega')$ is cohomologous to $(\phi,\omega)$ and $\psi\in\C^1(\frak g;\frak h)$ satisfies $\delta^1\psi=(\phi-\phi',\omega-\omega')$, then it is easy to see that that map $\alpha:(h,g)\mapsto (h+\psi(g),g)$ is an equivalence of restricted extensions.  These constructions are evidently inverse to each other so that we have shown the following theorem.
\begin{theorem}
\label{thm:rextalg}
The set of equivalence classes of restricted extensions of $\frak g$ by $\frak h$ is in one to one correspondence with $H^2(\frak g;\frak h)$.
\qed
\end{theorem}
An important special case is $\frak h=\mathbb F$.  
\begin{corollary}
$H^2(\frak g)$ is in one to one correspondence with equivalence classes of 1 dimensional restricted extensions of $\frak g$.
\qed
\end{corollary}
We remark here that the proof of the independence of the cohomology class on the equivalence class of the extension shows precisely why we require a strongly abelian algebra $\frak h$ in the definition of restricted extensions.  That is, it will always be the case that the skew-symmetric maps on $\frak g$ are (classically) cohomologous if they arise from different splitting maps.  However, the induced map $\tilde\psi$ need not equal the difference of the maps $\omega$ and $\omega'$ on $\frak g$ unless $\frak h$ is strongly abelian.  In terms of the canonical map $H^2(\frak g;\frak h)\to H^2_{\rm cl.}(\frak g;\frak h)$, we see that this map is injective precisely when $\frak h$ is strongly abelian.

We concluded this subsection with an investigation of the notion of a restricted infinitesimal deformation of a restricted Lie algebra. Since $\mathbb F$ has positive characteristic, our approach is purely algebraic.  

\begin{definition}
A restricted infinitesimal deformation of a restricted Lie algebra $\frak g$ is a skew-symmetric bilinear map $\phi:\frak g\times \frak g\to \frak g$ together with a map $\omega:\frak g\to \frak g$ such that for all $t\in\mathbb F$ the maps
\begin{equation}\label{defbracket}
\begin{array}{rcl}(g,h)&\mapsto &[gh]_t=[gh]+\phi(g,h)t \end{array}
\end{equation}
\begin{equation}\label{defpop}
\begin{array}{rcl}
g & \mapsto & g^{[p]_t}=g^{[p]}+\omega (g) t
\end{array}
\end{equation}
give the vector space $\frak g$ a restricted $\frak g$-module structure $\pmod {t^2}$.  
\end{definition}
Equivalently a restricted infinitesimal deformation is a restricted Lie algebra structure on the tensor product $(\mathbb F[t]/(t^2))\otimes_{\mathbb F} \frak g$ such that 
\[\epsilon\otimes 1_\frak g:(\mathbb F[t]/(t^2))\otimes \frak g\to \mathbb F \otimes \frak g=\frak g\]
is a restricted Lie algebra homomorphism where $\epsilon:\mathbb F[t]/(t^2)\to \mathbb F$ is the  canonical augmentation.  Two restricted infinitesimal deformations are equivalent if there is a linear map $\psi:\frak g\to \frak g$ such that 
\[\phi_1(g,h)-\phi_2(g,h)=[g\psi(h)]+[h\psi(g)]-\psi([gh])\]
and 
\[\omega_1(g)-\omega_2(g)=\psi(g^{[p]})-[\psi(g),\underbrace{g,\dots,g}_{p-1}].\]
It is well known that the bracket in (\ref {defbracket}) satisfies the Jacobi identity if and only if $\phi\in C^2_{\rm cl.}(\frak g;\frak g)$ is a cocycle. We claim that the properties of the $p$-operator in (\ref {defpop}) imply that $\omega$ has the $*$-property with respect to $\phi$ so that $(\phi,\omega)\in C^2(\frak g,\frak g)$ and that $(\phi,\omega)$ is a cocycle.  Indeed, since $(\lambda g)^{[p]_t}-\lambda^p g^{[p]_t}=0$ for all $t$, we easily see that $\omega(\lambda g)=\lambda^p \omega(g)$.  Since $t^2=0$, it is easy to expand the bracket $[g_1,\cdots,g_p]_t$ and we have
\[[g_1,\cdots,g_p]_t=[g_1,\cdots,g_p]+t\left (\sum_{k=0}^{p-2}[\phi([g_1,\dots,g_{p-k-1}],g_{p-k}),g_{p-k+1},\dots,g_p]\right ).
\]
Therefore, comparing the constant terms and the coefficients of $t$ in the identity
\[
(g+h)^{[p]_t}=g^{[p]_t}+h^{[p]_t}+ \sum_{\stackrel{g_j=g\ {\rm or}\  h}{\scriptscriptstyle g_1=g,g_2=h}} \frac{1}{\#(g)} [[[\cdots[[g_1g_2]_tg_3]_t\cdots]_tg_{p-1}]_tg_p]_t,
\]
we have the usual identity for $(g+h)^{[p]}$ and  
\begin{eqnarray*}
\lefteqn{\omega(g+h)=\omega(g)+\omega(h)}\\
& & +\sum_{\stackrel{g_j=g\ {\rm or}\  h}{\scriptscriptstyle g_1=g,g_2=h}} \frac{1}{\#(g)}\left (\sum_{k=0}^{p-2}[\phi([g_1,\dots,g_{p-k-1}],g_{p-k}),g_{p-k+1},\dots,g_p]\right ) .
\end{eqnarray*} 
Recalling that we have adjoint coefficients, we have 
\[\phi([g_1,\dots,g_{p-k-1}],g_{p-k}),g_{p-k+1},\dots,g_p]=(-1)^k g_p\cdots g_{p-k+1}\phi([g_1,\dots,g_{p-k-1}],g_{p-k})\]
so that $\omega$ has the $*$-property with respect to $\phi$ as claimed. It remains to show that the pair $(\phi,\omega)$ is a restricted cocycle.  Expanding the left hand side of the identity
\begin{equation}\label{conditionthree}
[gh^{[p]_t}]=[g,\underbrace{h,\dots,h}_p]_t,
\end{equation}
 and simplifying $\pmod {t^2}$, we have
\begin{equation}\label{lhs}
[gh^{[p]_t}]=[gh^{[p]}]+t\left ( [g\omega(h)]+\phi(g,h^{[p]})\right ).
\end{equation}
Expanding the right hand side of (\ref {conditionthree}), we have
\begin{eqnarray}
\lefteqn{[g,\underbrace{h,\dots,h}_p]_t=}\nonumber \\
& & [g,\underbrace{h,\dots,h}_p]+t\left (\sum_{i+j=p-1}[\phi([g,\underbrace{h,\dots,h}_j],h),\underbrace{h,\dots,h}_i]\right )\nonumber \\
& & [g,\underbrace{h,\dots,h}_p]+t\left (\sum_{i+j=p-1}(-1)^i h^i\phi([g,\underbrace{h,\dots,h}_j],h)\right ).\label{rhs}
\end{eqnarray}
Comparison of the coefficients of $t$ in (\ref {lhs}) and (\ref {rhs}) shows that 
\[\phi(g,h^{[p]})-\sum_{i+j=p-1}(-1)^i h^i\phi([g,\underbrace{h,\dots,h}_j],h)+[g\omega(h)]=0.\]
Therefore the induce map $\beta$ in (\ref {3induced}) is identically zero and the pair $(\phi,\omega)$ is a cocycle as claimed.  If two extensions $(\phi_1,\omega_1)$ and $(\phi_2,\omega_2)$ are equivalent and $\psi:\frak g\to \frak g$ is the map realizing this equivalence, then by definition $\delta^1\psi=(\phi_1-\phi_2,\omega_1-\omega_2)$ so that the resulting cohomology classes are the same.  Therefore we have a map from the set of equivalence classes of restricted infinitesimal deformations of $\frak g$ to $H^2(\frak g;\frak g)$. Conversely, if $(\phi,\omega)\in C^2(\frak g;\frak g)$ is a cocycle, and we define a a bracket $[,]_t$ and $p$-operator $\cdot^{[p]_t}$ on $\frak g$ by the formulae (\ref {defbracket}) and (\ref {defpop}), the the above arguments easily reverse to show that $\frak g$ is a restricted Lie algebra with these operations and hence we have a restricted infinitesimal deformation of $\frak g$.  Moreover, cohomologous cocycles give equivalent deformations by definition.  Summarizing, we have shown the following theorem.
\begin{theorem}
\label{thm:rinfinitesimaldeformations}
The equivalence classes of restricted infinitesimal deformations of a restricted Lie algebra $\frak g$ coincide with elements of $H^2(\frak g;\frak g)$.
\qed
\end{theorem}
We remark that our investigation of restricted infinitesimal deformations amounted to investigating the kernel of the canonical map $H^2(\frak g;\frak g)\to H^2_{\rm cl.}(\frak g;\frak g)$.  That is, two (ordinary) infinitesimal deformations of $\frak g$ in the sense of (2.1.4) may be equivalent via $\psi$, but $\tilde\psi$ need not satisfy (\ref {defpop}).

\section{Multiplicative Structures}
We conclude this chapter with a brief remark on the multiplicative structure of the complex $C$ defined in section (3.1).  

Recall that a differential graded algebra is a graded algebra $A$ over $\mathbb F$ equipped with degree $-1$ endomorphism $d:A\to A$ such that $d^2=0$ and the Leibniz formula
\[d(ab)=(da)b+(-1)^aa(db)\]
holds for all $a,b\in A$.  Here we use the symbol $a$ to denote both the element $a\in A$ as well as the degree of this element.  From this formula, we see immediately that the product of two cycles is again a cycle, and the product of a cycle $c$ and a boundary $db$ is a boundary $d(cb)$.  Therefore we can define the product of two homology classes as the class containing the product of any two representatives.  This product gives $H(A)$ itself the structure of a graded algebra.   Our purpose in bringing this up here is merely to mention that our complex $C$ is a differential graded algebra, and that this fact is useful in some of our computations above.  Namely, if we continue to fix a basis $\{e_1,\dots, e_n\}$ in $\frak g$, then as an algebra, $C_*=\bigoplus_{k\ge 0} C_k$ is generated by the elements
\begin{eqnarray*}
g^0_i & = & 1\otimes 1\otimes e_i \\
g^1_i & = & 1\otimes e_i\otimes 1 \\
g^2_i & = & e_i\otimes 1\otimes 1
\end{eqnarray*}
where $g^j_i\in C_j$ for $j=0,1,2$.  Using these generators, it is easy to check that the boundary map $d:C_*\to C_*$ satisfies the Leibniz formula
\[d(ab)=d(a)b+(-1)^a ad(b)\]
for all $a,b\in C_*$, and hence the homology $H(C)$ is a graded algebra.  In particular, we can now easily prove a claim made in section (3.1).  Recall that if we define 
\[c_i=1\otimes e_i^{[p]}\otimes 1-1\otimes e_i\otimes e_i^{p-1},\]
then $c_i\in C_1$ is a cycle for all $i=1,\dots,n$.  Since the map $d$ is multiplicative, it follows immediately that each product 
\[c_{i_1}\cdots c_{i_k}\]
is a cycle as claimed in section (3.1). We remark that without this multiplicative property of the boundary map $d$, the verification that $c_{i_1}\cdots c_{i_k}$ is a cycle is non-trivial.


\newpage
\pagestyle{myheadings} 
\markright{  \rm \normalsize CHAPTER 4. \hspace{0.5cm}
  COHOMOLOGY OF RESTRICTED LIE ALGEBRAS}
\chapter{Conclusions}
\thispagestyle{myheadings} 
In this concluding chapter, we briefly discuss some questions that remain unanswered by our work as well as give some indications of the methods we will pursue to answer these questions. 

First, in the case of an abelian restricted Lie algebra $\frak g$, we were able to construct a free augmented complex $C=\{C_k,d_k\}$
\[C_*\longrightarrow \mathbb F\longrightarrow 0\]
in the category $U_{\rm res.}(\frak g)$-modules that is exact in dimension less than $p$.  The cohomology of the derived complex $\Hom(C,M)$ therefore agrees with the Cartan-Eilenberg definition of cohomology in dimensions less than $p$.  We recall that in order to prove that this complex is exact for $k<p$, it was necessary to compute the homology of two auxiliary complexes $\C=\{\C_k,\partial_{\C}\}$ and $\frak C=\{\frak C_k,\partial_{\frak C}\}$.  We argued that given a cycle $c\in C_k$, we could modify $c$ by adding boundaries so that the leading term in $c$ was a cycle in the auxiliary complex $\frak C$.  This leading term is then homologous to a cycle of a special form so that we could write $c$ as a sum of a boundary and a cycle of lower total degree.  Induction on the degree of the leading term then showed that $c$ was in fact a boundary.  We recall here that the auxiliary boundary operators $\partial_{\C}$ and $\partial_{\frak C}$ are nothing more than the terms in the boundary operator $d$.  That is we have 
\[d=\partial_{\C}+\partial_{\frak C}.\]
All of this suggests that the auxiliary complexes $\C$ and $\frak C$ are nothing more than the initial terms of a spectral sequence that converges to the homology of the complex $C$.  That is, we believe that there is a filtration of the complex $C$ such that the initial terms of the corresponding spectral sequence are exactly the auxiliary complexes $\C$ and $\frak C$.  One goal in our subsequent investigations of this complex is to attempt to better organize the argumentation given in section (3.1) of this dissertation and find this filtration.  We remark that this by itself will not lead to any new results, as we were able to compute the homology of $C$ anyway.  However, the organization of such a filtration may very well lead to new insight into the structure of the complex $C$ and hence the resulting derived complex $\Hom(C,M)$.

As stated in section (3.2), our construction of the cochain complex 
\[ 0\to C^0({\frak g};M)\to C^1({\frak g};M)\to C^2({\frak g};M)\to C^3({\frak g};M)
\]
in the non-abelian case is sufficient for many computational purposes.  However, the argumentation involved in the constructions of these spaces and coboundary operators was extremely specialized and does not readily generalize to higher dimensions.  In our future research, we would like to obtain general constructions of the restricted cochain spaces $C^k(\frak g;M)$ and coboundary operators $\delta:C^k(\frak g;M)\to C^{k+1}(\frak g;M)$. As we remarked earlier, the equality of the dimensions of the cochain spaces constructed in the abelian and non-abelian cases (Corollary (\ref {cor:dimensions})) indicates that we may be able to develop the cohomology theory in the non-abelian case by deforming the spaces constructed in the abelian case. In fact, the deformation theory of restricted Lie algebras is a natural point to begin all of our cohomological investigations.

Recall that the connection between the cohomology of a Lie algebra $\frak g$ and infinitesimal deformations of $\frak g$ is encoded in $H^2_{\rm cl.}(\frak g;\frak g)$. A Lie algebra can be defined as a certain odd codifferential on the exterior coalgebra of a vector space, and by definition, an $L_\infty$ algebra is an arbitrary codifferential on this exterior coalgebra.  Using this point of view, it has been shown in \cite{P} that the ordinary Lie algebra cohomology $H_{\rm cl.}(\frak g;\frak g)$ classifies the infinitesimal deformations of the Lie algebra into an $L_\infty$ algebra.  Equivalently, we can say that the infinitesimal deformations of a Lie algebra into an $L_\infty$ algebra classifies $H_{\rm cl.}(\frak g;\frak g)$.  Conversations with the author of \cite{P} indicate that perhaps a restricted Lie algebra is merely a special case of some codifferential on the exterior coalgebra or the symmetric coalgebra (or a mixture of the two) of a vector space.  We have defined the restricted cohomology $H^2(\frak g;\frak g)$ of a restricted Lie algebra $\frak g$ and hence have defined the notion of infinitesimal deformations of these algebras.  It may therefore be possible to parallel the theory developed in \cite{P} for restricted Lie algebras.  We remark that it is often necessary to have $H^3(\frak g;\frak g)$ to complete the deformation theory and we do not have this space in the non-abelian case.  In fact, if the prime $p=3$, the complex in the abelian case fails to be exact at $k=3$ so that we do not have $H^3(\frak g;\frak g)$ is this case as well.  One of our primary goals in our future research will be to fully develop the deformation theory of restricted Lie algebras.  As we have already remarked, one consequence of this development may be to understand the notion of a restricted $L_\infty$ algebra.  Perhaps even more interesting is the fact stated in Corollary (\ref {cor:dimensions}): the dimensions of the cochain spaces in the abelian and non-abelian cases are identical.  This suggests that we may be able to construct the general non-abelian cochain spaces by deforming the spaces constructed in the abelian case.  In this sense the entire (non-abelian) restricted cohomology theory may be a quantization of the abelian restricted cohomology.  

Our future work will also include making computations of the cohomology of certain well known restricted Lie algebras.  In particular, we recall that if $\mathbb Z_p$ denotes the cyclic group of order $p$ and $A=\mathbb F(\mathbb Z_p)$ denotes the group algebra of $\mathbb Z_p$ over $\mathbb F$, then the Witt algebra is the derivation algebra $W=\Der A$. Recall that $W$ has a basis $D_j$, $(j=0,1,\dots, p-1)$, where 
\[D_jx=x^{j+1}\]
and  $W$ is a restricted Lie algebra with the operations
\begin{eqnarray*}
[D_i,D_j]&=&(j-i)D_{i+j} \\
D_0^{[p]}&=&D_0 \\
D_j^{[p]}&=&0\ (j>0).
\end{eqnarray*}
We remark that $W$ is simple as a Lie algebra so that the $p$-operator $D^{[p]}=D^p$ is the only map $W\to W$ giving $W$ a restricted Lie algebra structure.  The (ordinary) cohomology of the Witt algebra has been extensively studied since the early 1970's.  Although much is known about the cohomology spaces of this algebra, the results are usually very complicated. That is, $W$ has several exotic cohomology classes, and the existence of these classes usually depends on the prime $p$ being sufficiently large.  All of this suggests that we might benefit from examining the restricted cohomology of $W$.  The computation of the cohomology of $W$ is important since certain cohomology classes of this Lie algebra (in the characteristic zero case) correspond to characteristic classes of foliations of codimension 1 called Godbillon-Vey classes.  Of course, thanks to Chevalley, we can construct smooth Lie groups with given Lie algebras in positive characteristic as well and hence we may be able to develop the theory of Godbillon-Vey classes for these manifolds.

Finally, we remark that it is well known that the representation theories of Kac-Moody Lie algebras, quantum groups and restricted Lie algebras have striking similarities. Therefore, the construction of a suitably small cochain complex for restricted Lie algebra cohomology may simultaneously produce useful methods in the cohomology and representation theories of Kac-Moody algebras and quantum groups.  This is a very exciting research prospect.



\newpage
\pagestyle{myheadings} 
\markright{  \rm \normalsize BIBLIOGRAPHY. \hspace{0.5cm}
  COHOMOLOGY OF RESTRICTED LIE ALGEBRAS}

\end{document}